\documentclass[12pt,twoside]{amsart}

\usepackage{amsmath,amsfonts,amssymb,amscd,amsthm,amsxtra}

\usepackage[a4paper]{geometry}

\usepackage{enumerate,verbatim}
\usepackage{ifthen}

\usepackage{url,textcomp}
\newcommand{\no}{\textnumero}

\makeatletter
  \DeclareRobustCommand\em
         {\@nomath\em \ifdim \fontdimen\@ne\font >\z@
                       \upshape \else \slshape \fi}
\makeatother

\newtheoremstyle{mythm}
  {9pt}
  {9pt}
  {\slshape}
  {0pt}
  {\bfseries}
  {.}
  { }
  {\thmname{#1} \thmnumber{#2}\thmnote{ (#3)}}

\newtheoremstyle{mythmnonum}
  {9pt}
  {9pt}
  {\slshape}
  {0pt}
  {\bfseries}
  {.}
  { }
  {\thmname{#1}\thmnote{ (#3)}}

\theoremstyle{mythmnonum} 

\theoremstyle{mythm} 
\newtheorem{Theorem}{Theorem}[section]
\newtheorem{Lemma}[Theorem]{Lemma}

\newtheorem{Proposition}[Theorem]{Proposition}

\newtheorem{Corollary}[Theorem]{Corollary}

\theoremstyle{definition} 
\newtheorem{Definition}[Theorem]{Definition}
\newtheorem{RemDef}[Theorem]{Remark and Definition}

\newtheorem{Remark}[Theorem]{Remark}

\newtheorem{Example}[Theorem]{Example}
 
\newcommand{\file}[1]{{\url{#1}}}

\newcommand{\IC}{\mathbf{C}}                     
\newcommand{\IF}{\mathbf{F}}                     
\newcommand{\IN}{\mathbf{N}}                     
\newcommand{\IZ}{\mathbf{Z}}                     
\newcommand{\SG}{\mathfrak{S}}                   
\newcommand{\eps}{{\varepsilon}}                 
\renewcommand{\phi}{{\varphi}}                   
\newcommand{\ox}{\otimes}                        
\newcommand{\NCPartitions}{NC}                     
\newcommand{\NC}{\NCPartitions}                    
\newcommand{\NCI}{\NCPartitions I}                    
\newcommand{\abs}[1]{\left\lvert #1 \right\rvert}  
\newcommand{\bigabs}[1]{\bigl\lvert #1 \bigr\rvert}  
\newcommand{\biggabs}[1]{\biggl\lvert #1 \biggr\rvert}  

\newcommand{\norm}[1]{\left\lVert #1 \right\rVert} 
\newcommand{\bignorm}[1]{\bigl\lVert #1 \bigr\rVert}  
\newcommand{\biggnorm}[1]{\biggl\lVert #1 \biggr\rVert}  
\newcommand{\exchm}{\alg{E}} 
\newcommand{\exch}{\mbox{$\exchm$}} 
\DeclareMathOperator*{\freeprod}{\bigstar}    
\DeclareMathOperator*{\freeprodB}{{\bigstar_{\alg{B}}}}    
\DeclareMathOperator{\cn}{{\mathrm cn}}                     

\newlength{\tmpl}

\newcommand{\boxstar}{
\setlength{\tmpl}{\fboxsep}
\setlength{\fboxsep}{0pt}
\,\fbox{$\star$}\,
\setlength{\fboxsep}{\tmpl}
}

\newcommand{\Ord}{{\mathcal O}}                  
\newcommand{\bub}[1]{{\overset{\circ}{#1}}}      

\newcommand{\indep}{\perp\kern-3pt\perp}

\newcommand{\alg}[1]{{\mathcal{#1}}}                
\newcommand{\hilb}[1]{{\mathfrak{#1}}}                

\newcommand{\phit}{\tilde{\varphi}}
\newcommand{\psit}{\tilde{\psi}}
\newcommand{\phitt}{\tilde{\tilde{\varphi}}}

\numberwithin{equation}{section}
\newcommand{\paperchapter}{paper}

\allowdisplaybreaks{}

\newcommand{\WSC}{WSC}
\newcommand{\eWSC}{\mbox{$\widetilde{\text{WSC}}$}}

\begin{document}

\title[Cumulants in Noncommutative Probability Theory IV]
{Cumulants in Noncommutative Probability Theory IV. Noncrossing Cumulants:
  De~Finetti's Theorem and $L^p$-Inequalities}
\author{Franz Lehner}

\thanks{Supported by the Austrian Science Fund (FWF) Project \no{}R2-MAT}

\address{
Franz Lehner\\
In\-sti\-tut f{\"u}r Mathe\-ma\-tik C\\
Tech\-ni\-sche Uni\-ver\-si\-t{\"a}t Graz\\
Stey\-rer\-gas\-se 30, A-8010 Graz\\
Austria}
\email{lehner@finanz.math.tu-graz.ac.at}
\keywords{Cumulants, partition lattice, M{\"o}bius inversion, free probability,
  noncrossing partitions, noncommutative probability}

\subjclass{Primary  46L53, Secondary 05A18}

\date{\today}

\begin{abstract}
  De~Finetti's theorem states that any exchangeable sequence of classical random
  variables is conditionally i.i.d.\ with respect to some~$\sigma$-algebra. In this paper
  we prove a ``free'' noncommutative analog of this theorem, 
  namely we show that any noncrossing exchangeability system with a faithful state
  which satisfies a so called weak singleton condition 
  can be embedded into an free product with amalgamation
  over a certain subalgebra such that the interchangeable algebras remain
  interchangeable with respect to the operator-valued expectation. 
  Vanishing of crossing cumulants can be verified
  by checking a certain weak freeness condition and the weak singleton
  condition is satisfied e.g.\ when the state is tracial.
  The proof follows the classical proof of De~Finetti's theorem,
  the main technical tool being a noncommutative $L^p$-inequality for i.i.d.\ 
  sums of centered noncommutative random variables in noncrossing
  exchangeability systems.
\end{abstract}

\maketitle{}

\tableofcontents{}

De~Finetti's theorem~\cite{Kingman:1978:uses,ChowTeicher:1978:probability,DiaconisFreedman:1980:finite}
states that for any exchangeable sequence of random 
variables there exists a $\sigma$-algebra conditional on which the sequence
is i.i.d. 
There are various noncommutative versions of this theorem
\cite{Stormer:1969:symmetric,HudsonMoody:1975:locally,Hudson:1981:analogs,Petz:1990:definetti,AccardiLu:1993:continuous},
all of which involve tensor product constructions or other
commutativity conditions.
Indeed there is no hope to obtain a general De~Finetti's theorem
without imposing additional conditions.
In this \paperchapter{} we consider conditions under which a kind of
``most noncommutative'' version of De~Finetti's theorem holds, namely a 
characterization of exchangeability systems which can be written as an amalgamated
free product.

The only prerequisite for this \paperchapter{} is part~I of the series \cite{Lehner:2002:Cumulants1},
where exchangeability systems are introduced and many examples are discussed.
In Section~I.4.5 of that
\paperchapter{} we presented the amalgamated free product as an operator valued
exchangeability system with a conditional expectation~$\psi$. Composing this
conditional expectation with a state on the amalgamated subalgebra gives rise
to a scalar valued exchangeability system.

The question now is, under which conditions can an arbitrary exchangeability system
$\exch$ be written in this form?

An obvious necessary condition is that crossing cumulants must vanish,
because this is the case for the operator-valued amalgamated free cumulants and
the $\exch$-cumulants are simply the expectations of the latter, see
Section~I.3.6.
Another necessary condition is a certain \emph{weak singleton condition}.
The singleton condition introduced
in~\cite{BozejkoSpeicher:1996:interpolations} is too strong, because
together with the vanishing of crossing cumulants it actually implies freeness.
The weak singleton condition to be defined below however turns out to be the
right one and is automatically satisfied if the state is tracial.

The construction of the conditional expectation essentially follows the
classical proof,
namely by adjoining the algebra~$\alg{B}$ of permutation invariant
random variables to the initial algebra and extending to it the expectation
functional.
Moreover there it is possible to construct a conditional 
expectation $\psi_\infty$ onto $\alg{B}$ as the limit of symmetrizing maps
$\psi_N$.
There are certain technical issues regarding the faithfulness
of the extension in the non-tracial case. These are solved by a certain
Khinchin-type $L^p$-inequality which is of some independent interest.

Unfortunately we could not find an ``application'' of the characterization
obtained in this paper, except perhaps a new description of freeness
with amalgamation (called ``weak freeness''),
see Section~\ref{sec:GeneralCumulant4:WeakFreeness}.

The \paperchapter{} is organized as follows.

In Section~\ref{sec:GeneralCumulant4:Preliminaries} we collect a few
definitions and lemmas needed for the statement and the proof of the main
result.

In Section~\ref{sec:GeneralCumulant4:DeFinettiLemma}
we adapt a proof from \cite{AccardiLu:1993:continuous} to the
noncrossing situation. It shows that the conditional expectations $\psi_N$ evaluated
at words of the form $X_1^{(h(1))} X_2^{(h(2))}\dotsm X_n^{(h(n))}$
asymptotically factor according to the connected components of the kernel 
partition~$\ker h$, with an error term of order $1/N$, cf.\ the analogous
commutative result in \cite{DiaconisFreedman:1980:finite}).

In Section~\ref{sec:GeneralCumulants4:Lpinequality} we prove a strong law of
large numbers in noncrossing exchangeability systems 
under the assumption of the weak singleton condition.

In Section~\ref{sec:GeneralCumulant4:WeakFreeness} we discuss a certain weak
freeness condition and show that together with the weak singleton condition it
implies that crossing cumulants vanish and thus weak freeness is the same as
freeness with amalgamation.

\section{Preliminaries and statement of main result}
\label{sec:GeneralCumulant4:Preliminaries}

In this section we collect the necessary definitions and auxiliary results
needed later on. For details we refer to part I~\cite{Lehner:2002:Cumulants1}.

\subsection{Exchangeability Systems and Cumulants}

We recall first that a noncommutative probability space is a pair $(\alg{A},\phi)$
consisting of a  complex algebra $\alg{A}$ with unit~$I$ 
and a linear functional $\phi:\alg{A}\to\IC$  such that $\phi(I)=1$.  
An \emph{exchangeability system} $\exch=(\alg{U},\phit,\alg{J})$ for the
noncommutative probability space~$(\alg{A},\phi)$ consists of another
noncommutative probability space $(\alg{U},\phit)$ and an infinite
family $\alg{J}=(\iota_k)_{k\in\IN}$ of state-preserving embeddings
$\iota_k:\alg{A}\to \alg{A}_k\subseteq \alg{U}$, 
which we conveniently denote by $X\mapsto X^{(k)}$,
such that the image algebras $\alg{A}_j$ are \emph{interchangeable}
with respect to $\phit$:
for any family~$X_1,X_2,\dots,X_n\in \alg{A}$,
and for any choice of indices~$h(1),\dots,h(n)$
the expectation is invariant under any permutation~$\sigma\in\SG_\infty$
in the sense that
\begin{equation}
  \label{eq:GeneralCumulant4:phih=phisigmah}
  \phit(X_1^{(h(1))} X_2^{(h(2))} \cdots X_n^{(h(n))})
  = \phit(X_1^{(\sigma(h(1)))} X_2^{(\sigma(h(2)))} \cdots X_n^{(\sigma(h(n)))})
  .
\end{equation}
Denote by~$\Pi_n$ the lattice of partitions (or equivalence relations) of the
set $[n]=\{1,2,\dots,n\}$. The refinement order $\pi\leq\sigma$ means as usual
that the partition $\pi$ is finer than the partition $\sigma$.
Then permutation invariance means that the value of the
expectation~\eqref{eq:GeneralCumulant4:phih=phisigmah} 
only depends on the so-called \emph{kernel} of the map~$h$
which is the partition~$\pi=\ker h\in\Pi_n$ defined by 
$$
i\sim_\pi j
\qquad
\iff
\qquad
h(i)=h(j)
$$
and we denote the common value as
\begin{equation}
  \label{eq:GeneralCumulant4:PartitionedExpectation}
  \phi_\pi(X_1,X_2,\dots,X_n) = 
  \phit(X_1^{(\pi(1))} X_2^{(\pi(2))} \cdots X_n^{(\pi(n))})
  .
\end{equation}
Here we consider a partition~$\pi\in\Pi_n$ as a function $\pi:[n]\to\IN$,
mapping each element to the number of the block containing it.
This is a canonical example of an index function~$h$ with $\ker h=\pi$
and because of condition~\eqref{eq:GeneralCumulant4:phih=phisigmah} the
actual numbering of the blocks does not matter.

Throughout this \paperchapter{} all algebras will be $C^*$- or
pre-$C^*$-algebras and we will assume that the algebra $\alg{U}$ is generated
by the algebras $\alg{A}_k$ and that the action of $\SG_\infty$ extends to
all of $\alg{U}$ leaving the state $\phit$ invariant.
For an index set $I\subseteq \IN$ we denote by $\alg{A}_I$
the algebra generated by $\{\alg{A}_i,i\in I\}$. While the state $\phi$ on $\alg{A}$ 
usually will be assumed to be faithful, this is not always true for the state
$\phit$ on $\alg{U}$. Indeed a major part of this paper is dedicated to
the proof that a certain GNS-state is at least partially faithful.

The constructions above can be done in the more general situation of an
\emph{operator-valued noncommutative probability space}, which is a pair
$(\alg{A},\psi)$ consisting of a unital algebra $\alg{A}$ and a
\emph{conditional expectation} $\psi$ onto some unital subalgebra
$\alg{B}\subseteq\alg{A}$.
Here a conditional expectation is a unital positive map $\psi:\alg{A}\to\alg{B}$,
with the property that $\psi(BXB')=B\psi(X)B'$ whenever $B,B'\in \alg{B}$ and
$X\in\alg{A}$. The free amalgamated exchangeability system is an example
of this more general concept, see below.

\begin{Example}
The most commutative example of an exchangeability system for an arbitrary noncommutative
probability space $(\alg{A},\phi)$ is the infinite tensor product
$$
\alg{U}=\bigotimes_{i=1}^\infty \alg{A}_i
$$
of infinitely many copies  $\alg{A}_i$ of $\alg{A}$ with the tensor product
state $\phit  = \otimes_{i=1}^\infty \phi_i$,
where $(\alg{A}_i,\phi_i)_{i=1}^\infty$ is an infinite family of copies
of $(\alg{A},\phi)$ and the embeddings are
$$
\iota_j:X\mapsto X^{(j)} = I\ox I\ox \cdots \ox I \ox X \ox I\ox\cdots
$$
Then the subalgebras are clearly interchangeable and the partitioned expectation
\eqref{eq:GeneralCumulant4:PartitionedExpectation}
evaluates to
$$
\phi_\pi(X_1,X_2,\dots,X_n) = \prod_{B\in\pi} \phi(\prod_{i\in B} X_i)
$$
which is familiar from classical probability theory.
\end{Example}
\begin{Example}
  Taking the reduced free product instead of the tensor product leads to the
  free exchangeability system
  $$
  (\alg{U},\phit) = \freeprod_{i=1}^\infty (\alg{A}_i,\phi_i)
  $$
  More generally, if $(\alg{A},\phi)$ in addition comes with a conditional
  expectation $\psi$ onto some subalgebra $\alg{B}$ such that $\phi\circ\psi=\phi$,
  then one can construct the amalgamated free
  exchangeability system
  $$
  (\alg{U},\psit) = \freeprodB_{i=1}^\infty{}
  (\alg{A}_i,\psi_i)  
  $$
  While $(\alg{U},\psit,\alg{J})$ is an operator valued exchangeability system
  for the operator valued noncommutativity space $(\alg{A},\psi)$,
  it becomes a scalar exchangeability system for $(\alg{A},\phi\circ\psi)$
  for any state $\phi$ on $\alg{B}$ by letting $\phit=\phi\circ\psit$.
\end{Example}
More examples are listed in \cite{Lehner:2002:Cumulants1}.
In some sense (made precise in \cite{Speicher:1997:universal}) the reduced free
product and the tensor product together with Boolean independence 
are the only universal exchangeability systems. The emphasis in \cite{Speicher:1997:universal}
however lies on ``universality'' in the sense that the partitioned moment functionals
$\phi_\pi(X_1,X_2,\dots,X_n)$ depend on the individual distributions of the $X_i$ 
in a universal way (i.e., as a polynomial formula). This already excludes the
free amalgamated exchangeability system constructed above;
our approach is less constructive as we assume
that an exchangeability system is given \emph{apriori} and 
we do not assume universality. The concept of ``identical distribution'' becomes
more involved, as explained below.

Subalgebras $\alg{B},\alg{C}\subseteq \alg{A}$ are called
\emph{$\alg{E}$-exchangeable} or, more suggestively, \emph{$\alg{E}$-independent}
if for any choice of random variables~$X_1,X_2,\dots,X_n \in \alg{B}\cup\alg{C}$
and subsets~$I,J\subseteq \{1,\dots,n\}$ such that~$I\cap J=\emptyset$,
$I\cup J=\{1,\dots,n\}$, $X_i\in \alg{B}$ for $i\in I$  and $X_i\in \alg{C}$ for $i\in J$,
we have the identity
$$
\phi_\pi(X_1,X_2,\dots,X_n)
= \phi_{\pi'}(X_1,X_2,\dots,X_n)
$$
whenever~$\pi$, $\pi'\in\Pi_n$ are partitions with~$\pi|_I=\pi'|_I$ and
$\pi|_J=\pi'|_J$.
Two families of random variables~$(X_i)$ and~$(Y_j)$ are called
$\alg{E}$-exchangeable if the algebras they generate have this
property.

We say that two random variables~$X$ and~$Y\in\alg{A}$
\emph{have the same distribution given~$\exchm$},
if for any word $W=W_1W_2\cdots W_n$ with
$W_i\in \{X^{(1)}\}\cup \bigcup_{i\geq2} \alg{A}_i$ the expectation~$\phit(W)$
does not change if we replace each occurrence of~$X^{(1)}$ by~$Y^{(1)}$.
We call~$X$ and~$Y$ \emph{\exch-i.i.d.}\ if in addition they are \exch-independent.
Similarly a sequence $(X_i)_{i\in\IN}\subseteq\alg{A}$ of \exch-independent random variables
is called \emph{\exch-i.i.d.}\
if for any word $W=W_1W_2\cdots W_n$ with
$W_i\in \{X_i^{(1)}: i\in\IN\}\cup \bigcup_{i\geq2} \alg{A}_i$ the expectation~$\phit(W)$
does not change if we apply a permutation $\sigma\in\SG_\infty$ to the indices
of $X_i$, i.e., if we replace each occurrence of $X_i$ by $X_{\sigma(i)}$.

Then it is possible to define cumulant functionals, indexed by set partitions
$\pi\in\Pi_n$, via
$$
K_\pi^\exchm(X_1,X_2,\dots,X_n)
= \sum_{\sigma\leq\pi}
   \phi_\sigma(X_1,X_2,\dots,X_n)\,\mu(\sigma,\pi) 
$$
where $\mu(\sigma,\pi)$ is the M{\"o}bius function of the lattice of set
partitions, cf.\ part~I.
The use of the probabilistic terminology ``independence'' and ``cumulants'' is
justified by the following proposition which establishes the analogy to
classical probability.

\begin{Proposition}[{\cite{Lehner:2002:Cumulants1}}]
  Two subalgebras~$\alg{B},\alg{C}\subseteq\alg{A}$ are \exch-independent 
  if and only if mixed cumulants vanish, that is,
  whenever~$X_i\in\alg{B}\cup\alg{C}$ are some noncommutative random variables
  and~$\pi\in\Pi_n$ is an arbitrary partition
  such that there is a block of $\pi$ which contains indices $i$ and $j$
  such that~$X_i\in\alg{B}$ and~$X_j\in\alg{C}$,
  then $K^\exchm_\pi(X_1,X_2,\dots,X_n)$ vanishes.
\end{Proposition}

With this abstract formalism one can transfer many combinatorial proofs from
classical probability to the general situation. One of the most useful
results is the product formula of Leonov and Shiryaev.

\begin{Proposition}[{\cite[Prop.~3.3]{Lehner:2002:Cumulants1}}]
  \label{prop:GeneralCumulants:Productformula}
  Let~$(X_{i,j})_{i\in\{1,\dots,m\},j\in \{1,\dots,n_i\}}\subseteq\alg{A}$
  be a family of noncommutative random variables
  containing in total $n=n_1+n_2+\dots+n_m$ variables.
  Then every partition~$\pi\in\Pi_m$ induces a partition~$\tilde{\pi}$
  on $\{1,\dots,n\}\simeq\{(i,j):i\in [m], j\in[n_i]\}$ with blocks
  $\tilde{B}=\{(i,j):i\in B, j\in[n_i]\}$, that is, each block~$B\in\pi$ is
  replaced by the union of the intervals~$(\{n_{i-1}+1,n_{i-1}+2,\dots,n_i\})_{i\in B}$.
  Then we have
  $$
  K^\exchm_\pi(\prod_{j_1} X_{1,j_1},\prod_{j_2} X_{2,j_2},\dots,\prod_{j_m} X_{m,j_m})
  =\sum_{\substack{\sigma\in\Pi_n\\ \sigma\lor \tilde{\hat0}_m=\tilde\pi}}
    K^\exchm_\sigma(X_{1,1},X_{1,2},\dots,X_{m,n_m})
  $$
\end{Proposition}

\begin{RemDef}
  \label{rem:GeneralCumulant4:iid=interchangeable}
  In the sequel we will frequently appeal to the following simple observation
  in order to reduce the amount of indices, see e.g.\   
  Corollary~\ref{cor:GeneralCumulant4:exchCSdifferentindices}.
  We will be dealing with noncommutative polynomials involving variables
  $X_1^{(h(1))}X_2^{(h(2))}\dotsm X_n^{(h(n))}$ for which we want to produce
  ``independent'' copies, that is, replacing $X_j^{(h(j))}$ by $X_j^{(h'(j))}$ in
  such a way that the ranges of the indices $h$ and $h'$ are disjoint.
  This can be interpreted as follows.
  Let $I$ be an index set containing all the indices $h(j)$,
  $j\in\{1,\dots,n\}$. Consider the algebra $\tilde{\alg{A}}=\alg{A}_I$
  generated by $(\alg{A}_i)_{i\in I}$ and a sequence $(I_j)_{j\in\IN}$ of
  mutually disjoint index sets $I_j\subseteq \IN$ of the same cardinality as $I$.
  Then the \emph{extended exchangeability system}
  $\tilde{\exchm}=(\alg{U},\phit,\tilde{\alg{J}})$ is an 
  exchangeability system for $\tilde{\alg{A}}$, where 
  $\tilde{\iota}_j:\tilde{\alg{A}}\to \tilde{\alg{A}}_j=\alg{A}_{I_j}$ is the natural
  permutation isomorphism induced by an arbitrary bijection between $I$ and $I_j$.
  The procedure of choosing $X_j^{(h'(j))}$ and $X_j^{(h''(j))}$ in the
  exchangeability system $\exch$ such that the ranges of $h'$ and
  $h''$ are disjoint amounts to the same as taking interchangeable copies $\tilde{X}_j^{(1)}$
  and $\tilde{X}_j^{(2)}$ of $\tilde{X}_j=X_j^{(h(j))}\in\tilde{A}$ 
  in the exchangeability system $\tilde{\exch}$.
  We will denote these by $X^{(I_j)}=\tilde{X}^{(j)}$.
  Let us illustrate this idea by a simplified example on the free group
  $\IF_\infty$. The group algebra of $\IF_\infty$ is an exchangeability
  system for the group algebra of $\IZ$ and at the same time
  it is an exchangeability system for the group algebra of $\IF_N$ for
  arbitrary $N$, because it can be written as
  $$
  \IF_\infty = \IZ*\IZ*\cdots =   \IF_N* \IF_N*\cdots
  .
  $$

  Thus we will sometimes do proofs for the initial exchangeability system
  $\exchm$ and state the results for $\tilde{\exchm}$ as corollaries.
  Also cumulants of polynomials $W_j\in \alg{A}_I$ are defined in $\tilde{\exchm}$
  as well; the values do not depend on the choice of the index sets $I$ and $I_j$.

  Similarly we will sometimes not distinguish between i.i.d.\
  sequences in~$\alg{A}$ in the sense defined earlier in this section 
  and sequences of the form~$X^{(i)}$.

\end{RemDef}

\subsection{Noncrossing partitions and freeness}

The lattice of noncrossing partitions, denoted by $\NC_n$, will play a
prominent r{\^o}le in this \paperchapter{}. We recall that a partition~$\pi\in\Pi_n$ is
\emph{noncrossing}, if there is no quadruple of indices $i<j<k<l$, such that
$i\sim_\pi k$ and $j\sim_\pi l$ and $i\not\sim j$.
Equivalently, noncrossing partitions can also be characterized recursively
by the property that there is always at least one block which is an interval
and after removing this block the remaining partition is still noncrossing.

\begin{Definition}
We say that \emph{crossing cumulants vanish} in a given exchangeability
system~$\exch$
if for any $n$ and for any choice of random variables $X_1$,
$X_2$,\ldots,~$X_n$ we have the identity
$$
K_\pi^\exchm(X_1,X_2,\dots,X_n) = 0
$$
whenever $\pi$ has a crossing, i.e., $\pi\not\in\NC_n$.
We call such an exchangeability system a \emph{noncrossing exchangeability system}.
\end{Definition}
A prominent example of a noncrossing exchangeability system 
is the free exchangeability system, where~$(\alg{U},\phit)$ 
is the reduced free product of an infinite family of
copies of a given noncommutative probability space $(\alg{A},\phi)$,
see Section~I.4.4.
We recall that the free exchangeability system is characterized by the
property that 
$$
\phit(X_1^{(h(1))} X_2^{(h(2))} \dotsm X_n^{(h(n))}) = 0
$$
whenever $\phi(X_j)=0$ and $h(j)\neq h(j+1)$ for every $1\leq j\leq n-1$.

Another situation where crossing cumulants vanish is \emph{freeness with amalgamation}
\cite{Voiculescu:1995:operations,Speicher:1998:combinatorial},
which is a noncommutative analog of conditional independence.
Let $(\alg{A},\psi)$ be a $\alg{B}$-valued noncommutative probability space.
Then the amalgamated free product
$(\alg{U},\psit)=\bigstar_{\alg{B}}(\alg{A}_i,\psi_i)$
of infinitely many copies of
$\alg{A}$ with amalgamation over $\alg{B}$ is a $\alg{B}$-valued noncrossing
exchangeability system for $(\alg{A},\psi)$.
The amalgamated free exchangeability system is characterized by the property
that
$$
\psit(X_1^{(h(1))} X_2^{(h(2))} \dotsm X_n^{(h(n))}) = 0
$$
whenever $\psi(X_j)=0$ and $h(j)\neq h(j+1)$ for every $1\leq j\leq n-1$.
As in the case of scalar freeness, the~$\alg{B}$-valued cumulants
$
K_\pi^{\exchm,\psi}(X_1,\dots,X_n)
$ 
vanish for any partition~$\pi\not\in\NC_n$.

Now choose any state $\phi$ on~$\alg{B}$, then
\begin{equation}
  \label{eq:GeneralCumulant4:exchph=U,phiopsi}
  \exch^\phi=(\alg{U}=\bigstar_\alg{B}\alg{A}_i,\phi\circ\psit,\alg{J})
\end{equation}
is a scalar-valued
exchangeability system for the noncommutative probability space
$(\alg{A},\phi\circ\psi)$ whose cumulants are
$$
K^{\exchm,\phi\circ\psi}_\pi(X_1,\dots,X_n)
=\phi(K^{\exchm,\psi}_\pi(X_1,\dots,X_n))
;
$$
again cumulants vanish for any partition~$\pi$ with crossings.
Other examples of noncrossing exchangeability systems can be constructed
by taking conditionally free
products~\cite{BozejkoLeinertSpeicher:1996:convolution}.
We will see later that these cannot be realized as amalgamated free products,
cf.~Remark~\ref{rem:GeneralCumulant3:ConditionalFreeness}

\subsection{Multiplicative functions and convolution on the lattice of
  noncrossing partitions}
We refer to \cite{Stanley:1986:Enumerative1} or
Section~I.1.3 for 
the definition of the incidence algebra of a poset.
In the case of noncrossing partitions there is also a reduced incidence algebra
of multiplicative functions.
Let~$\NCI_n$ be the set of intervals in $\NC_n$ and $\NCI=\bigcup_n \NCI_n$.
Then every interval $[\pi,\sigma]$ has a canonical decomposition
\begin{equation}
  \label{eq:GeneralCumulant4:pisigma=decomposition}
[\pi,\sigma] \simeq [\hat0_1,\hat1_1]^{k_1} \times[\hat0_2,\hat1_2]^{k_2}
\times\dotsm\times[\hat0_1,\hat1_n]^{k_n} \times \dotsm
\end{equation}
where~$(k_n)_{n=1}^\infty$ is a sequence of integers with finitely many nozero
entries and~\eqref{eq:GeneralCumulant4:pisigma=decomposition} is a lattice
isomorphism~\cite{Speicher:1994:multiplicative}.
For example it is easy to see that for~$\pi\in\NC_n$ we have
\begin{equation}
  \label{eq:GeneralCumulant4:0pi=decomposition}
  [\hat0_n,\pi] \simeq   \prod_{p=1}^n [\hat0_p,\hat1_p]^{k_p}
\end{equation}
where~$k_p$ is the number of blocks of~$\pi$ of size~$p$. 

A function $f:\NCI\to\IC$ is called \emph{multiplicative} if for any interval
$[\pi,\sigma]$
it satisfies
$$
f([\pi,\sigma]) 
= f([\hat0_1,\hat1_1])^{k_1} f([\hat0_2,\hat1_2])^{k_2} \dotsm f([\hat0_1,\hat1_n])^{k_n}
$$
where $[\pi,\sigma]$ has the
decomposition~\eqref{eq:GeneralCumulant4:pisigma=decomposition}.
Such a function is determined by its \emph{characteristic sequence}~$f_n=f([\hat0_n,\hat1_n])$
and it is easy to see that the convolution of two multiplicative functions is
again multiplicative, i.e., the multiplicative functions constitute an algebra,
the so-called \emph{reduced incidence algebra}.
For a noncrossing partition~$\pi$ with decomposition as
in~\eqref{eq:GeneralCumulant4:0pi=decomposition} we will denote 
$$
f_\pi := f([\hat0_n,\pi])= \prod_{p=1}^n f_p^{k_p}
.
$$
Then the convolution~$f\boxstar g$ of two multiplicative functions can be
calculated with the aid of the \emph{Kreweras complementation map}
\cite{Kreweras:1972:partitions}.
This is a lattice anti-automorphism of~$\NC_n$ described as follows.
Paint $n$ points on a circle and label them clockwise with numbers
$1$,~$2$,~\ldots,~$n$. A noncrossing partition $\pi\in\NC_n$ can be visualized
by drawing inside the circle for each block of $\pi$ the convex polygon whose vertices are the
elements of the block. Now put another~$n$ points with labels
$\bar1$,~$\bar2$,~\ldots,~$\bar{n}$ on the circle, placing the point with label
$\bar{k}$ between the points with label $k$ and $k+1$ and connect the new
points with each other by drawing as many lines as possible without
intersecting the polygons drawn before. This leads to a noncrossing partition
of the set $\{\bar1,\bar2,\dots,\bar{n}\}$ which is called the 
\emph{Kreweras complement} of~$\pi$ and is denoted~$K(\pi)$.
It is easy to see that~$K$ is an order anti-automorphism, $K(\hat0_n)=\hat1_n$ and
$K(\hat1_n)=\hat0_n$.
It follows that $[\pi,\hat1_n]\simeq [\hat0_n,K(\pi)]$ and therefore
the convolution of multiplicative functions can be written
$$
(f\boxstar g)_n=f\boxstar g([\hat0_n,\hat1_n])
=\sum_{\pi\in\NC_n} f_\pi g_{K(\pi)}
.
$$
As a consequence the reduced incidence algebra is commutative and there is a
``Fourier Transform'' \cite{NicaSpeicher:1997:fourier}:
Let $(f_n)_{n=1}^\infty$, $(g_n)_{n=1}^\infty$ be the characteristic sequences
of two multiplicative functions $f$ and $g$ with $f_1=g_1=1$.
The formal power series
$$
\phi_f(z) = \sum_{n=1}^\infty f_n z^n
$$
is called the \emph{characteristic series} of $f$.
Let $\phi_f^{\langle-1\rangle}(z)$ be the compositional inverse of $\phi_f$ and
$$
\alg{F}_f(z) = \frac{1}{z} \, \phi_f^{\langle-1\rangle}(z),
$$
then 
\begin{equation}
  \label{eq:GeneralCumulant4:Ffboxstarg=FfFg}
  \alg{F}_{f\boxstar g}(z) = \alg{F}_f(z)\,\alg{F}_g(z)
\end{equation}
Two prominent multiplicative functions are the \emph{Zeta function}
$\zeta(\pi,\sigma)\equiv 1$ with characteristic series 
$$
\phi_\zeta(z)=\frac{z}{1-z}
$$
and its inverse, the \emph{M{\"o}bius function} $\mu$ whose
characteristic sequence is given by the signed Catalan numbers
$\mu_n=(-1)^{n-1}C_{n-1}=\frac{(-1)^{n-1}}{n}\binom{2n-2}{n-1}$ and the
characteristic series is
$$
\phi_{\mu}(z) = \frac{\sqrt{1+4z}-1}{2}
.
$$
These functions satisfy $\zeta\boxstar\mu=\delta$, where $\delta$ is the unit
element of the reduced incidence algebra and has characteristic series
$\phi_\delta(z)=z$.
Moreover, functions $f$ and $g$ satisfy $f\boxstar \zeta= g$ if and only if
$f=g\boxstar \mu$ and this is the case if and only if
\begin{equation}
  \label{eq:GeneralCumulant4:fboxstarzeta=g}
  \phi_f(z(1+\phi_g(z))) = \phi_g(z)
  .  
\end{equation}
We will encounter applications of these formulae in
Section~\ref{sec:GeneralCumulants4:Lpinequality}.

\subsection{The weak singleton condition}

One more ingredient is needed for the formulation of the main result.
We have already seen that the vanishing of crossing cumulants is a necessary
condition for an exchangeability system to come from an amalgamated free product.
This condition however is not sufficient as will be shown below,
namely a so called \emph{weak singleton condition} is also necessary.

\begin{Definition}[{\cite{BozejkoSpeicher:1996:interpolations}}]
  \begin{enumerate}[(a)]
   \item 
  An exchangeability system $\exch=(\alg{U},\phit,\alg{J})$ satisfies the 
  \emph{singleton condition} if
  $$
  \phit(X_1^{(i_1)} X_2^{(i_2)} \dotsm X_n^{(i_n)}) = 0
  $$
  whenever one of the indices $i_j$ appears only once and the
  corresponding random variable $X_j$ satisfies $\phi(X_j)=0$.
 \item   An exchangeability system $\exch$ satisfies the 
  \emph{weak singleton condition} (\WSC{}) if
  $$
  \phit(X_1^{(i_1)} X_2^{(i_2)} \dotsm X_n^{(i_n)}) = 0
  $$
  whenever one of the indices $i_j$ appears only once and the
  corresponding random variable $X_j$ satisfies $\phit({X_j^{(1)}}^*X_j^{(2)})=0$.
 \item   An exchangeability system $\exch$ satisfies the 
  \emph{extended weak singleton condition} (\eWSC{}) if
  the extended exchangeability system $\tilde{\exchm}$ of
  Definition~\ref{rem:GeneralCumulant4:iid=interchangeable} satisfies (\WSC{}), i.e.,
  for any finite index set $I\subseteq\IN$ and any $n$-tuple of polynomials
  $W_1,W_2,\dots,W_n\in \alg{A}_I$ we have
  $$
  \phit(W_1^{(I_{i_1})} W_2^{(I_{i_2})} \dotsm W_n^{(I_{i_n})}) = 0
  $$
  whenever one of the index sets $I_{i_j}$ appears only once and the
  corresponding polynomial $W_j$ satisfies $\phit({W_j^{(I_1)}}^*W_j^{(I_2)})=0$.
  \end{enumerate}
\end{Definition}

\begin{Remark}
  The weak singleton condition is indeed weaker than the singleton condition,
  because it follows from
  Corollary~\ref{cor:GeneralCumulant4:ExchangeableCauchySchwarz}
  that the condition~$\phit({X^{(1)}}^*X^{(2)})=0$ implies $\phi(X)=0$.

  The extended \WSC{} is introduced for technical reasons and needed only for
  Corollary~\ref{cor:GeneralCumulants4:KhinchinGeneral}.
  We were not able to prove or disprove that it follows from (\WSC{}) in general,
  however it is automatically implied by (\WSC{}) if the initial algebra
  contains ``enough'' independent random variables, i.e., if for
  any $n$-tuple $(X_1,\dots,X_n)$ of elements of~$\alg{A}$
  there exist arbitrary many i.i.d.\ copies inside~$\alg{A}$.
  This is the case in all examples known to the author.
\end{Remark}

The next proposition shows that a weak singleton condition holds also for cumulants,
if it holds for moments.
\begin{Proposition}
  \label{prop:GeneralCumulant4:WeakSingletonCumulants}
  Let $\exch$ be an exchangeability system in which the weak singleton
  condition holds.
  Then $K_\pi^{\exchm}(X_1,X_2,\dots,X_n)=0$ whenever the partition~$\pi$
  contains a singleton $\{j\}$ such that $\phit({X_j^{(1)}}^*X_j^{(2)})=0$.
\end{Proposition}
\begin{proof}
  Indeed,
  $$
  K_\pi^{\exchm}(X_1,X_2,\dots,X_n)
  =\sum_{\sigma\leq\pi} \phi_\sigma(X_1,X_2,\dots,X_n)
   \,\mu(\sigma,\pi)
  $$
  and all terms vanish, because each $\sigma\leq\pi$ contains the singleton $\{j\}$.
\end{proof}
The starting point of this \paperchapter{} is the following observation and its
corollary.
\begin{Lemma}
  \label{lem:GeneralCumulants4:singletonpartition}
  Let $\pi\in\Pi_n$ be an alternating partition, i.e., a partition in which
  neighbouring elements are in different blocks. Then any noncrossing partition $\sigma\leq\pi$
  contains at least one singleton.
\end{Lemma}
\begin{proof}
  Any noncrossing partition $\sigma$ contains at least one interval block and
  the condition $\sigma\leq\pi$ implies that this interval block has
  length~$1$, i.e., it is a singleton.
\end{proof}

\begin{Corollary}
  \label{cor:GeneralCumulant4:singleton+crossing=free}
  A noncrossing exchangeability system which satisfies the singleton condition
  is given by a reduced free product. 
\end{Corollary}

\begin{proof}
  Let $X_1,\dots,X_n\in\alg{A}$ with $\phi(X_j)=0$ and let
  $h(1)$, $h(2)$, \ldots, $h(n)$ be indices such that $h(j)\neq h({j+1})$.
  We have to show that $\phit(X_1^{(h(1))} X_2^{(h(2))}\dotsm X_n^{(h(n))})=0$.
  The singleton condition implies that $\phi_\pi(X_1,X_2,\dots,X_n)=0$
  whenever $\pi$ contains a singleton and consequently for any such $\pi$ the
  corresponding cumulant~$K_\pi^{\exchm}(X_1,X_2,\dots,X_n)$ vanishes.
  Now 
  $$
  \phit(X_1^{(h(1))} X_2^{(h(2))}\dotsm X_n^{(h(n))})
  =\sum_{\pi\leq\ker h}
    K_\pi^{\exchm}(X_1,X_2,\dots,X_n)
  $$
  and by assumption the sum extends over noncrossing partitions only.
  By Lemma~\ref{lem:GeneralCumulants4:singletonpartition} any such
  partition contains a singleton and the corresponding cumulant vanishes
  because of Proposition~\ref{prop:GeneralCumulant4:WeakSingletonCumulants}.
\end{proof}

\subsection{Conditional Expectations}
\label{subsec:GeneralCumulant4:ConditionalExpectations}

The proof of~\ref{cor:GeneralCumulant4:singleton+crossing=free}
stays essentially the same if the singleton condition is replaced by (\eWSC{}),
resulting in an amalgamated free product.
The main technical problem is the construction of the conditional expectation~$\psi$ onto a
certain algebra~$\alg{B}$ and to prove faithfulness of an extension of~$\phit$
on~$\alg{B}$ in order to apply Lemma~\ref{lem:GeneralCumulants4:singletonpartition}.
The construction of $\psi$ is the same as in the commutative case, namely
as the limit of symmetrizing maps.
\begin{Definition}
  \label{def:GeneralCumulant4:psiN}
  Let $\exch=(\alg{U},\phit,J)$ be an exchangeability system for some
  noncommutative probability space $(\alg{A},\phi)$.
  We define the conditional expectations $\psi_N$, $N\in\IN$, by
  $$
  \psi_N(X)=\frac{1}{N!} \sum_{\sigma\in\SG_N} \sigma(X)
  ;
  $$
  for polynomials, i.e., elements of the form $X=X_1^{(h(1))} X_2^{(h(2))}\dotsm
  X_n^{(h(n))}$ this is
  $$
  \psi_N(X_1^{(h(1))} X_2^{(h(2))}\dotsm X_n^{(h(n))})
  = \frac{1}{N!}
    \sum_{\sigma\in\SG_N}
     X_1^{(\sigma(h(1)))} X_2^{(\sigma(h(2)))}\dotsm X_n^{(\sigma(h(n)))}
  $$
  Note that if $h$ is fixed and $N$ large enough, this only depends on $\ker h$.
\end{Definition}
We collect a few elementary properties of $\psi_N$. Proofs are easy and can be found
in~\cite{AccardiLu:1993:continuous}.
\begin{Proposition}
  \begin{enumerate}
   \item $\phit\circ\psi_N=\phit$.
   \item $\psi_N(X)=X$ if and only if $\sigma(X)=X$ $\forall\sigma\in\SG_N$.
   \item $\psi_N\circ\psi_M=\psi_N$ for $M\leq N$.
   \item $\psi_N\circ \iota_k = \frac{1}{N}\sum_{j=1}^N \iota_j$ if $k\leq N$.
  \end{enumerate}
\end{Proposition}

In contrast to finite exchangeability systems, there is a nonnegative bilinear
form available in the infinite case.
\begin{Proposition}
  The sesquilinear form
  $$
  \langle X,Y\rangle=\phit({Y^{(1)}}^* X^{(2)})
  $$
  is nonnegative on~$\alg{A}$.
\end{Proposition}
\begin{proof}
  Consider for fixed $N\in\IN$ the nonnegative expectation
  \begin{align*}
    \phit(\psi_N(X^{(1)})^*\psi_N(X^{(1)}))
    &= \left(
         \frac{1}{N!}
       \right)^2
       \sum_{\sigma,\sigma'} 
        \phit({X^{(\sigma(1))}}^* X^{(\sigma'(1))}) \\
    &= \frac{1}{N^2}
       \sum_{i,j=1}^N \phit({X^{(i)}}^*{X}^{(j)})\\
    &= \frac{N(N-1)}{N^2}\,\phit({X^{(1)}}^*X^{(2)}) + \frac{1}{N}
    \phit({X^{(1)}}^*X^{(1)})
    .
  \end{align*}
  Now letting $N\to\infty$ yields the claim.
\end{proof}
Positivity implies the Cauchy-Schwarz inequality.
\begin{Corollary}
  \label{cor:GeneralCumulant4:ExchangeableCauchySchwarz}
  For any $X,Y\in\alg{A}$
  $$
  \bigabs{\phit({Y^{(1)}}^* X^{(2)})}
  \leq \phit({X^{(1)}}^* X^{(2)})^{1/2}\,\phit({Y^{(1)}}^* Y^{(2)})^{1/2}
  $$
\end{Corollary}
Similarly one can prove a multivariable Cauchy-Schwarz inequality,
cf.\ Remark~\ref{rem:GeneralCumulant4:iid=interchangeable}.
\begin{Corollary}
  \label{cor:GeneralCumulant4:exchCSdifferentindices}
  Let $\{g(1),g(2),\dots,g(m)\}$ and $\{h(1),h(2),\dots,h(n)\}$ be
  disjoint sets of indices and
  $P=P(X_1^{(g(1))},X_2^{(g(2))},\dots,X_m^{(g(m))})$ and
  $Q=P(Y_1^{(h(1))},Y_2^{(h(2))},\dots,Y_m^{(h(n))})$
  noncommutative polynomials in $X_i$ and $Y_i$,
  then 
  $$
  \abs{\phit(Q^*P)}^2 \leq \phit(P^*P')\,\phit(Q^*Q')
  $$
  where 
  $
  P'=P(X_1^{(g'(1))},X_2^{(g'(2))},\dots,X_m^{(g'(m))})
  $ 
  and 
  $
  Q'=P(Y_1^{(h'(1))},Y_2^{(h'(2))},\dots,Y_m^{(h'(n))})
  $ 
  such that
  $h'$ (resp.~$g'$) is an index function whose range is disjoint
  from the range of $h$ (resp.~$g$).
\end{Corollary}
After these preparations we can consider two situations in which the weak 
singleton condition holds.
\begin{Proposition}
  \label{prop:GeneralCumulant4:weaksingletoncondition}
  Each of the following two conditions implies (\WSC) (and (\eWSC)).
  \begin{enumerate}[(a)]
   \item The state is faithful and the exchangeability system comes from an amalgamated free
    product as described in   \eqref{eq:GeneralCumulant4:exchph=U,phiopsi}.
   \item The state is tracial.
  \end{enumerate}
\end{Proposition}
\begin{proof}
  \begin{enumerate}[(a)]
   \item Let $\psit$ be the conditional expectation with respect to
    which the algebras $\alg{A}_i$ are free. Then by
    Proposition~\ref{prop:Generalcumulant4:Speicherspsipi} below we have
    \begin{align*}
      0 &= \phit({X^{(1)}}^* X^{(2)}) \\
        &= \phit(\psit({X^{(1)}}^* X^{(2)})) \\
        &= \phit(\psit(X^{(1)})^*\psit(X^{(1)}))
    \end{align*}
    and by faithfulness this implies that $\psit(X^{(1)})=0$
    Now if $X$ appears as a singleton in some
    word, then the expectation of the word vanishes. Indeed, if~$X_j=X$ and the
    index $h(j)$ appears only once in the range of the index function $h$,
    then we may condition on $\alg{A}_{h(j)}$ (see
    Proposition~\ref{def:Generalcumulant4:Speicherspsipi} below) and obtain
    $$
    \phit(X_1^{(h(1))} X_2^{(h(2))}\dotsm X_n^{(h(n))})
    = \phit(X_1^{(h(1))} X_2^{(h(2))}\dotsm\psit(X_j^{(h(j))})\dotsm X_n^{(h(n))})
    = 0
    $$
   \item Let $X_1$, $X_2$,\ldots,$X_n\in\alg{A}$ and
    $h:[n]\to\IN$ be an index function such that $h(j)$ is a
    singleton and assume that~$\phit({X_j^{(1)}}^*X_j^{(2)})=0$. 
    We have to show that $\phit(X_1^{(h(1))}\dotsm X_n^{(h(n))})$ vanishes.
    By traciality we may
    assume without loss of generality that $j=n$.
    Then we may apply the Cauchy-Schwarz inequality of
    Corollary~\ref{cor:GeneralCumulant4:exchCSdifferentindices} and
    with any index function $h'$ whose range is disjoint from that of
    $h$ we obtain
    \begin{align*}
    \bigabs{
      \phit(X_1^{(h(1))}\dotsm X_n^{(h(n))})
    }^2
    &\leq \phit(X_1^{(h(1))}
               \dotsm
               X_{n-1}^{(h(n-1))}
               {X_{n-1}^{(h'(n-1))}}^*
               \dotsm
               {X_{1}^{(h'(1))}}^* )
          \,
          \phit({X_n^{(h'(n))}}^* X_n^{(h(n))})\\
    &=0
    \end{align*}
  \end{enumerate}
\end{proof}

\subsection{Statement of main result}

The first part of
Proposition~\ref{prop:GeneralCumulant4:weaksingletoncondition} shows
that the weak singleton condition is a necessary condition for an
exchangeability system to come from an amalgamated free product.
We can now state the main theorem of this \paperchapter{}.
\begin{Theorem}
  \label{thm:GeneralCumulant4:MainTheorem}
  Let $(\alg{A},\phi)$ be a noncommutative probability space and 
  $\exch=(\alg{U},\phit,\alg{J})$ a noncrossing exchangeability system
  for $(\alg{A},\phi)$ with faithful state $\phit$ which satisfies
  (\eWSC{}). Then $\exch$ 
  can be embedded into a $\alg{B}$-valued exchangeability system 
  $\tilde{\alg{E}}=(\tilde{U},\alg{J},\psi)$ such that the interchangeable algebras $\alg{A}_i$
  are free with amalgamation over $\alg{B}$ and interchangeable with respect to $\psi$.
\end{Theorem}
\begin{Remark}
  While it is true that any exchangeability system can be embedded into an
  amalgamated free product (the trivial one, where $\alg{B}$ coincides with
  the full algebra), it is not always true that this can be done in such a way
  that the $\alg{A}_i$ are still interchangeable. Therefore the preceding
  theorem is nontrivial. This is like in the commutative case,
  where an arbitrary exchangeable sequence of random variables is
  trivially conditionally independent with respect to the full $\sigma$-algebra, 
  but they are certainly not conditional i.i.d, unless they are identical.
\end{Remark}
\begin{Remark}
  \label{rem:GeneralCumulant3:ConditionalFreeness}
  Other examples where crossing cumulants vanish are Boolean independence
  \cite{SpeicherWoroudi:1997:boolean}
  and more generally conditional free
  independence~\cite{BozejkoLeinertSpeicher:1996:convolution}.
  In these examples however Theorem~\ref{thm:GeneralCumulant4:MainTheorem} does
  not apply because either the state is not faithful
  or the weak singleton condition fails:
  Let~$(\alg{U},\phit,\psi)=\bigstar(\alg{A}_i,\phi_i,\psi_i)$ be the
  conditionally free exchangeability system for $(\alg{A},\phi)$,
  cf.~\cite{BozejkoLeinertSpeicher:1996:convolution} or
  Section~I.4.7.
  For our purposes it is sufficient to know the defining property
  $$
  \phit(X_1^{(h(1))}X_2^{(h(2))}\dotsm X_n^{(h(n))})=
  \phi(X_1)\,\phi(X_2)\dotsm \phi(X_n)
  $$
  whenever $\psi(X_j)=0$ and $h(j)\neq h(j+1)$ for every $1\leq j\leq n-1$,
  and the modified pyramidal law
  $$
  \phit(X_1^{(1)}Y^{(2)}X_2^{(1)})
  = \phi(X_1)\,\phi(Y)\,\phi(X_2) + \psi(Y)\,(\phi(X_1X_2) - \phi(X_1)\,\phi(X_2))
  ,
  $$
  cf.~\cite{BozejkoLeinertSpeicher:1996:convolution} or
  Lemma~I.4.14.
  Let us investigate the weak singleton condition for some
  element~$X\in\alg{A}$.
  Denoting~$\bub{X}=X-\psi(X)$ we compute
  \begin{align*}
    \phit({X^{(1)}}^*X^{(2)})
    &= \phit( (\bub{X}{{}^{(1)}}^* + \psi(X^*) )\,(\bub{X}{{}^{(2)}} + \psi(X)))\\
    &=(\phit( \bub{X}{{}^{(1)}}^*\bub{X}{{}^{(2)}}) +
    (\phi(X^*)-\psi(X^*))\,\psi(X)+ \psi(X^*)\,(\phi(X)-\psi(X)) +
    \abs{\psi(X)}^2\\
    &= \abs{\phi(X)}^2
  \end{align*}
  Now let $X\in\alg{A}$ be any element with $\phi(X)=0$ but $\psi(X)\neq0$
  (this is possible unless~$\phi=\psi$; in the latter case we have just usual
  freeness)
  and find elements $Y$ and $Z$ such that $\phi(YZ)\neq\phi(Y)\,\phi(Z)$.
  Then
  $$
  \phit(Y^{(1)}X^{(2)}Z^{(1)}) =
  \phi(Y)\,\phi(X)\,\phi(Z) + \psi(X)\,(\phi(YZ)-\phi(Y)\,\phi(Z))
  $$
  does not vanish as it should if the weak singleton condition were true.
  It follows from Proposition~\ref{prop:GeneralCumulant4:weaksingletoncondition}
  that the conditional free product cannot be embedded into a
  free amalgamated exchangeability system with a faithful state~$\phit$.
\end{Remark}

In the remaining sections we will construct a conditional expectation $\psi$
on $\alg{U}$ and show that the algebras $\alg{A}_i$ are free with respect to this $\psi$.
The latter is constructed by a law of large numbers, namely
as limit of the symmetrizing maps $\psi_N$ of Definition~\ref{def:GeneralCumulant4:psiN}.
This is motivated by the following heuristics.
Assume that $\alg{A}_i$ are free with respect to some
some conditional expectation $\psi$, then it is known that for any
$X\in\alg{A}$ with $\psi(X)=0$ the norm
\begin{equation}
  \label{eq:Generalcumulant4:amalgamatednorminequality}
  \biggnorm {\sum_{i=1}^N X^{(i)} } \leq 2\sqrt{N}\norm{X}
\end{equation}
and therefore
$$
\frac{1}{N}\sum X^{(i)} = \psi(X) + \frac{1}{N}\sum (X^{(i)}-\psi(X))
$$
converges to $\psi(X)$ in norm as $N$ tends to infinity.
We will prove an inequality similar to
\eqref{eq:Generalcumulant4:amalgamatednorminequality}
in Section~\ref{sec:GeneralCumulants4:Lpinequality} by combinatorial methods,
i.e., without assuming freeness and using only (\WSC{}) and the fact that crossing cumulants
vanish.

\section{A De~Finetti Lemma}
\label{sec:GeneralCumulant4:DeFinettiLemma}

In this section we prove an asymptotic factorization property of the
conditional expectations of Definition~\ref{def:GeneralCumulant4:psiN}.
First we need to review noncrossing partitioned conditional
expectations~\cite{Speicher:1998:combinatorial}.
\begin{Definition}[{\cite{Speicher:1998:combinatorial}}]
\label{def:Generalcumulant4:Speicherspsipi}
  Let~$\psi$ be a conditional expectation.
  For a noncrossing partition $\pi\in\NC_n$ let~$b=\{k,k+1,\dots,l\}$
  be an interval block and define recursively
  $$
  \psi[\pi](X_1,X_2,\dots,X_n)
  = \psi[\pi\textstyle{\setminus} b](X_1,
                                 X_2,
                                 \dots,
                                 X_{k-1},
                                 \psi(X_k X_{k+1} \dotsm X_l) X_{l+1},
                                 \dots,
                                 X_n)
  $$
\end{Definition}
These partitioned expectations appear in the calculation of Speicher's amalgamated free
cumulants.
\begin{Proposition}[{\cite{Speicher:1998:combinatorial}}]
  \label{prop:Generalcumulant4:Speicherspsipi}
  Let~$\exch=(\bigstar_\alg{B}\alg{A}_i,\tilde{\psi},\alg{J})$
  be the amalgamated free exchangeability system for a $\alg{B}$-valued
  noncommutative probability space $(\alg{A},\psi)$.
  Then for any noncrossing partition~$\rho$ 
  and any finite sequence
  $X_1$,~$X_2$,\ldots,~$X_n\in\alg{A}$ the partitioned expectations
  \eqref{eq:GeneralCumulant4:PartitionedExpectation} coincide with Speicher's
  partitioned expectations in Definition~\ref{def:Generalcumulant4:Speicherspsipi}:
  $$
  \psi_\rho(X_1,X_2,\dots,X_n)
  = \psi[\rho](X_1,X_2,\dots,X_n)
  .
  $$
  More generally, for an arbitrary index function~$h$ and any noncrossing
  partition~$\rho$ s.t.\ $\rho\geq\ker h$ we have
  $$
  \tilde{\psi}(X_1^{(h(1))} X_2^{(h(2))}\dotsm X_n^{(h(n))})
  =  \tilde{\psi}[\rho](X_1^{(h(1))} X_2^{(h(2))}\dotsm X_n^{(h(n))})
  $$
\end{Proposition}
We will show that the conditional expectations of
Definition~\ref{def:GeneralCumulant4:psiN} asymptotically have the same property.
For the proof of this fact the following elementary estimate is needed in two places.
\begin{Lemma}
  \label{lem:GeneralCumulants4:1m1mjdN}
  Let $j$, $p$, $N$ be positive integers with $j\leq N$ and $p\leq N$,
  then
  \begin{equation}
    \label{eq:GeneralCumulants4:1m1mjdN}    
    1-
    \left(
      1-\frac{j}{N}
    \right)
    \left(
      1-\frac{j}{N-1}
    \right)
    \dotsm
    \left(
      1-\frac{j}{N-p+1}
    \right)
    \leq\frac{pj}{N-p+1}
  \end{equation}
\end{Lemma}
\begin{proof}
  Denote $E_p$ the left hand side of \eqref{eq:GeneralCumulants4:1m1mjdN}.
  Clearly the sequence $E_p$ satisfies
  $0\leq E_p\leq1$, is nondecreasing and therefore~$E_p\geq E_1=\frac{j}{N}$.
  Moreover it satisfies the recursion
  \begin{align*}
    E_{p+1}
    &= 1-
       \left(
         1-\frac{j}{N}
       \right)
       \left(
         1-\frac{j}{N-1}
       \right)
       \dotsm
       \left(
         1-\frac{j}{N-p}
       \right)\\
    &= E_p + \frac{j}{N-p}(1-E_p)
  \end{align*}
  We proceed by induction to show that $C_p=\frac{pj}{N-p+1}$ is an upper bound.
  Suppose that for $E_p$ we the estimate $E_p\leq C_p$ holds.
  Then
  \begin{align*}
    E_{p+1}
    &\leq C_p+\frac{j}{N-p}
       \left(
         1-\frac{j}{N}
       \right)\\
    &= \frac{pj}{N-p+1}
       +
       \frac{j}{N-p}
       \left(
         1-\frac{j}{N}      
       \right)\\
    &\leq \frac{pj}{N-p} + \frac{j}{N-p}\\
    &\leq \frac{(p+1)j}{N-p}
  \end{align*}
\end{proof}

The proof of the following inequality has been adapted to noncrossing
partitions from~\cite[Lemma~2.6]{AccardiLu:1993:continuous};
the estimate goes back to and is a noncommutative analog
of the main result in~\cite{DiaconisFreedman:1980:finite}.

\begin{Lemma}
  \label{lem:GeneralCumulants4:DeFinettiLemma}
  Let $\pi\in\Pi_n$ and $\rho\in\NC_n$ s.t.\ $\rho\geq\pi$
  containing $p=\abs{\pi}$ and $r=\abs{\rho}$ blocks, respectively.
  Then for $N\geq p$
  $$
  \norm{\psi_N(X_1^{(\pi(1))} X_2^{(\pi(2))}\dotsm X_n^{(\pi(n))})
    -\psi_N[\rho](X_1^{(\pi(1))}, X_2^{(\pi(2))},\dotsm, X_n^{(\pi(n))})}
  \leq\frac{(2r-1)p^2}{N-p+1} \prod\norm{X_i}
  $$
\end{Lemma}
\begin{proof}
  Let~$\rho=\rho_1<\rho_2<\dots<\rho_r=\hat{1}_n$ be a maximal chain in
  $[\rho,\hat{1}_n]\cap\NC_n$ with the property that the blocks $b_j$ of $\rho$ can 
  be labeled in such a way that
  $$
  \rho_k=\{b_1\cup b_2\cup\dots\cup b_k,b_{k+1},\dots,b_r\}
  .
  $$
  Such a chain can be constructed by ordering the chains with respect to their
  minimal elements and then successively merging the leftmost two blocks.
  Correspondingly we label the blocks $a_j$ of~$\pi$ in such a way that
  $a_1,\dots,a_{j_1}\subseteq b_1$, 
  $a_{j_1+1},\dots,a_{j_2}\subseteq b_2$, etc.,
  $a_{j_{r-1}+1},\dots,a_{j_r}\subseteq b_r$.
  Denote~$p_k=j_k-j_{k-1}$ the number of blocks of~$\pi$ which are contained
  in the~$k$th block~$b_k$ of~$\rho$.
  Let $\tilde{\rho}=\rho/\pi$, i.e., the partition of the block set of
  $\pi$ induced by $\rho$: For $i$, $j\in\{1,\dots,p\}$ we set
  $i\sim_{\tilde{\rho}}j$ if $a_i$ and $a_j$ are contained in the same
  block of $\rho$. 
  We have to compare the first term
  \begin{equation}
    \label{eq:GeneralCumulant4:accardilufirstexpression}
    \begin{aligned}
      \psi_N(X_1^{(\pi(1))} X_2^{(\pi(2))}\dotsm X_n^{(\pi(n))})
      &=\frac{1}{N!}
        \sum_{\sigma\in\SG_N}
         X_1^{(\sigma(\pi(1)))}
         X_2^{(\sigma(\pi(2)))}
         \dotsm
         X_n^{(\sigma(\pi(n)))}\\
      &=\frac{(N-p)!}{N!}
        \sum_{\substack{h:[p]\to[N]\\ \ker h =\hat{0}_p}}
         X_1^{(h(\pi(1)))}
         X_2^{(h(\pi(2)))}
         \dotsm
         X_n^{(h(\pi(n)))}
    \end{aligned}
  \end{equation}
  with the second one
  \begin{equation}
    \label{eq:GeneralCumulants4:definettilemmasecondterm}
    \begin{aligned}
    \psi_N[\rho](X_1^{(\pi(1))} X_2^{(\pi(2))}\dotsm X_n^{(\pi(n))})
    &=\left(
        \frac{1}{N!}
      \right)^r
      \sum_{\sigma_1,\dots,\sigma_r\in\SG_N}
       X_1^{(\sigma_{\rho(1)}(\pi(1)))}
       X_2^{(\sigma_{\rho(2)}(\pi(2)))}
       \dotsm
       X_n^{(\sigma_{\rho(n)}(\pi(n)))}\\
    &=\biggl(
        \prod_{k=1}^r \frac{(N-p_k)!}{N!}
      \biggr)
      \sum_{\substack{h:[p]\to[N]\\ \ker h\land\tilde{\rho} =\hat{0}_p}}
       X_1^{(h(\pi(1)))}
       X_2^{(h(\pi(2)))}
       \dotsm
       X_n^{(h(\pi(n)))}
    \end{aligned}
  \end{equation}
  As in the proof of~\cite[Lemma~2.6]{AccardiLu:1993:continuous} we now
  split~\eqref{eq:GeneralCumulants4:definettilemmasecondterm} as
  \begin{multline*}
    \psi_N[\rho](X_1^{(\pi(1))} X_2^{(\pi(2))}\dotsm X_n^{(\pi(n))})\\
    =\biggl(
        \prod_{k=1}^r \frac{(N-p_k)!}{N!}
      \biggr)
      \biggl(
        \sum_{k=1}^{r-1}
         \sum_{\substack{h:[p]\to[N]\\ \ker h\land\tilde{\rho}_k =\hat{0}_p\\ \ker h\land\tilde{\rho}_{k+1} >\hat{0}_p}}
          X_1^{(h(\pi(1)))}
          X_2^{(h(\pi(2)))}
          \dotsm
          X_n^{(h(\pi(n)))}
\\
         +
         \sum_{\substack{h:[p]\to[N]\\ \ker h =\hat{0}_p}}
          X_1^{(h(\pi(1)))}
          X_2^{(h(\pi(2)))}
          \dotsm
          X_n^{(h(\pi(n)))}
       \biggr)
  \end{multline*}
  Up to a multiplicative constant, the last term is the same
  as~\eqref{eq:GeneralCumulant4:accardilufirstexpression}, and we will
  show that the constants are asymptotically the same;
  but first we will bound the remaining $r-1$ terms of the sum.
  The conditions $\ker h\land\tilde{\rho}_k=\hat{0}_p$ and  $\ker
  h\land\tilde{\rho}_{k+1}>\hat{0}_p$
  mean that $h|_{\{1,\dots,j_k\}}$ is injective, but $h|_{\{1,\dots,j_{k+1}\}}$ is
  not,
  i.e., at least one of the indices $h(j_k+1),h(j_k+2),\dots,h(j_{k+1})$ is
  contained in $\{h(1),\dots,h(j_k)\}$; remember that
  $h|_{\{j_k+1,\dots,j_{k+1}\}}$ \emph{is} injective.
  Thus
  $$
  \sum_{\substack{h:[p]\to[N]\\
                  \ker h\land\tilde{\rho}_k =\hat{0}_p\\ 
                  \ker h\land\tilde{\rho}_{k+1} >\hat{0}_p}}
  =\sum_{h(1),\dots,h(k) \text{ distinct}}
   \sum_{\substack{h(j_k+1),\dots,h(j_{k+1}) \text{ distinct}\\ 
       \{h(1),\dots,h(k)\}\cap \{h(j_k+1),\dots,h(j_{k+1})\}\neq\emptyset}}
   \sum_{h(j_{k+2}),\dots,h(j_r)}
  $$
  There are $N(N-1)\dotsm(N-j_k+1)$ different choices for~$h(1)$,~$h(2)$,\ldots,$h(j_k)$,
  $$
  N(N-1)\dotsm(N-p_{k+1}+1)-(N-j_k)(N-j_k-1)\dotsm(N-j_k-p_{k+1}+1)
  $$
  choices for~$h(j_k+1)$,~$h(2)$,\ldots,$h(j_{k+1})$, and
  $$
  \prod_{s=k+2}^r \frac{N!}{(N-p_s)!}
  $$
  possibilities to choose the remaining indices~$h(j_{k+1}+1)$,\ldots,$h(p)$.
  The $k$th term can therefore be estimated by
  $$
  \biggl(
    \prod_{s=1}^{k+1} \frac{(N-p_s)!}{N!}
  \biggr)
  \frac{N!}{(N-j_k)!}
  \left(
    \frac{N!}{(N-p_{k+1})!} - \frac{(N-j_k)!}{(N-j_k-p_{k+1})!}
  \right)
  \prod_{j=1}^n\norm{X_j}
  $$
  By Lemma~\ref{lem:GeneralCumulants4:1m1mjdN}
  \begin{align*}
    \frac{(N-p_{k+1})!}{N!}
    &
    \left(
      \frac{N!}{(N-p_{k+1})!}
      -
      \frac{(N-j_k)!}{(N-j_k-p_{k+1})!}
    \right)\\
    &= 1-\frac{(N-j_k)(N-j_k-1)\dotsm(N-j_k-p_{k+1}+1)}{N(N-1)\dotsm(N-p_{k+1}+1)}\\
    &= 1-
    \left(
      1-\frac{N-j_k}{N} \frac{N-j_k-1}{N-1}\dotsm\frac{N-j_k-p_{k+1}+1}{N-p_{k+1}+1}
    \right)\\
    &= 1-
    \left(
      1-\frac{j_k}{N}
    \right)
    \left(
      1-\frac{j_k}{N-1}
    \right)
    \dotsm
    \left(
      1-\frac{j_k}{N-p_{k+1}+1}
    \right)\\
    &\leq \frac{p_{k+1}j_k}{N-p_{k+1}+1}
  \end{align*}
  and therefore the $k$th term is smaller than
  \begin{multline*}
    \frac{N(N-1)\dotsm(N-j_k+1)}%
         {N(N-1)\dotsm(N-p_1+1)\dotsm N(N-1)\dotsm(N-p_k+1)}
    \frac{p_{k+1}j_k}{N-p_{k+1}+1}
    \prod\norm{X_i}\\
    \leq
    \frac{p_{k+1}j_k}{N-p_{k+1}+1}
    \prod\norm{X_i}
  \end{multline*}
  and
  $$
  \sum_{k=1}^{r-1}
  \frac{p_{k+1}j_k}{N-p_{k+1}+1}
  \leq(r-1)\frac{\bar{p}p}{N-\bar{p}+1}
  $$
  where $\bar{p}=\max p_k$.
  Now we come to the difference between the final term
  and~\eqref{eq:GeneralCumulant4:accardilufirstexpression}.
  \begin{align*}
    \bignorm{
      \biggl(
        \prod_{k=1}^r \frac{(N-p_k)!}{N!}
        &
        -\frac{(N-p)!}{N!}
      \biggr)
      \sum_{\ker h=\hat{0}_p }
       X_1^{(h(\pi(1)))}
       X_2^{(h(\pi(2)))}
       \dotsm
       X_n^{(h(\pi(n)))} 
    }\\
    &\leq 
     \biggl(
       \prod_{k=1}^r \frac{(N-p_k)!}{N!}
       -\frac{(N-p)!}{N!}
     \biggr)
     N(N-1)\dotsm(N-p+1)
     \prod\norm{X_i}\\
    &=
     \biggl(
       1 - \prod_{k=1}^r \frac{(N-p_k)!}{N!} N(N-1)\dotsm(N-p+1)
     \biggr)
     \prod\norm{X_i}\\
    &=
     \biggl(
       1 - \prod_{k=1}^r \frac{(N-j_{k-1})(N-j_{k-1}-1)\dotsm(N-j_k+1)}{N(N-1)\dotsm(N-p_k+1)}
     \biggr)
     \prod\norm{X_i}\\
    &=
     \biggl(
       1
       -
       \prod_{k=1}^r
        \left(
          1-\frac{j_{k-1}}{N}
        \right)
        \left(
          1-\frac{j_{k-1}}{N-1}
        \right)
        \dotsm
        \left(
          1-\frac{j_{k-1}}{N-p_k+1}
        \right)
     \biggr)
     \prod\norm{X_i}\\
    &\leq
     \biggl(
       1
       -
       \prod_{k=1}^r
        \left(
          1-\frac{p_k j_{k-1}}{N-p_k+1}
        \right)
     \biggr)
     \prod\norm{X_i}\\
    &\leq
    \left(
      1- \left(
            1-\frac{p\bar{p}}{N-\bar{p}+1}
          \right)^r
    \right) \prod\norm{X_i}\\
    &\leq
      r \frac{p\bar{p}}{N-\bar{p}+1}
      \prod\norm{X_i}
  \end{align*}
  by Lemma~\ref{lem:GeneralCumulants4:1m1mjdN}.
\end{proof}

The conditional expectations $\psi_N$ need not converge but we can construct a
limit by extending the algebra with the help of the GNS-construction as
in~\cite{AccardiLu:1993:continuous}; the price of this is a possible loss
of faithfulness, which will be repaired in the next section.
Let $\pi:\alg{U}\to B(\hilb{H})$ be the GNS representation of $\alg{U}$ on 
$\hilb{H}=L^2(\alg{U},\phit)$. By assumption it is faithful and cyclic with
cyclic vector $\xi_0$, i.e.,
$\{\pi(X)\,\xi_0 : X\in\alg{U}\}$ is a dense subspace of~$\hilb{H}$
and $\phit(X)=\langle \pi(X)\,\xi_0,\xi_0\rangle$.
Since we assumed that $\alg{U}$ is generated by $(\alg{A}_i)_{i\in I}$,
the action of $\SG_\infty$ on $\alg{U}$ can be extended to a representation~$U_\sigma$
on~$\hilb{H}$ which is characterized by
$$
U_\sigma\pi(X)\,\xi_0 = \pi(\sigma(X))\,\xi_0
\qquad\forall\sigma\in\SG_\infty\forall X\in\alg{U}
$$
Let 
$$
\hilb{H}_\infty = \{\xi \in\hilb{H} : U_\sigma\xi = \xi\quad\forall\sigma\in\SG_\infty\}
$$
be the subspace of $U$-invariant elements and $P_\infty:\hilb{H}\to\hilb{H}_\infty$ the
orthogonal projection, and define
$$
\psi_\infty(X) = P_\infty\pi(X)P_\infty
.
$$
Similarly the projection~$P_N$ onto 
$$
[\pi(\psi_N(\alg{U}))\,\xi_0]=\{\xi\in\hilb{H}:U_\sigma\xi=\xi
\quad\forall\sigma\in\SG_N\}
$$ 
is characterized by the property
$$
P_N\pi(X)\,\xi = \pi(\psi_N(X))\,\xi
\qquad\forall\sigma\in\SG_N\forall X\in\alg{U}
$$
Clearly every $P_N\geq P_\infty$ and the sequence $P_N$ is monotonically
decreasing to its strong limit $P_\infty$. 
We continue our work in the extended noncommutative probability
space~$(\tilde{\alg{U}},\phitt)$ generated by 
$\pi(\alg{U})$ and $P_\infty$ and where the state $\phitt(X) = \langle
X\xi_0,\xi_0\rangle$
is the GNS-extension of $\phit$.
We may also  consider it as an operator-valued
noncommutative probability space with the conditional expectation
$$
\psi_\infty:\tilde{\alg{U}}\to \alg{B}=P_\infty \tilde{\alg{U}}P_\infty
$$
and we have~$\phitt=\phitt\circ\psi_\infty$.
Moreover,
$$
\phitt(\psi_\infty(X_1)\psi_\infty(X_2)\dotsm\psi_\infty(X_n))
= \lim_{N\to\infty}
   \phit(\psi_N(X_1)\psi_N(X_2)\dotsm\psi_N(X_n))
;
$$
in particular, $\phitt$ is a trace if $\phit$ is a trace.
As a corollary to Lemma~\ref{lem:GeneralCumulants4:DeFinettiLemma} we have the
following generalization of~\cite[Lemma~3.1]{AccardiLu:1993:continuous}. 

\begin{Lemma}
  \label{lem:GeneralCumulants4:DeFinettiLemmaGNS}
  Let $h:[n]\to\IN$ be an index function and
  let $\rho$ be any noncrossing partition such that $\rho\geq\ker h$.
  Then for any sequence $X_1$,~$X_1$,\ldots,~$X_n\in\alg{A}$ we have the factorization
  \begin{align*}
  \phit(X_1^{(h(1))}X_2^{(h(2))}\dotsm X_n^{(h(n))})
  &= \lim_{N\to\infty}
     \langle 
       \psi_N[\rho](X_1^{(h(1))},X_2^{(h(2))},\dots,X_n^{(h(n))})\,\xi_0,\xi_0
     \rangle\\
  &= \langle 
       \psi_\infty[\rho](X_1^{(h(1))},X_2^{(h(2))},\dots,X_n^{(h(n))})\,\xi_0,\xi_0
     \rangle
  \end{align*}
\end{Lemma}
\begin{proof}%
  [First part of the proof of Theorem~\ref{thm:GeneralCumulant4:MainTheorem}]
  Let $\alg{U}_0\subseteq\tilde{\alg{U}}$ be the algebra of polynomials,
  i.e., the (non-closed) algebra generated by $(\alg{A}_i)_{i\in I}$ and
  $\alg{B}_0= \psi_\infty(\alg{U}_0)$ its image under $\psi$ (as a vector space).
  Let us assume for a moment that~$\phitt$ is faithful on~$\alg{B}_0$
  in the sense that for any element $W \in \alg{B}_0$ the equation
  $\phitt(W^*W)=0$ implies that $W=0$.
  We show that the images~$\tilde{\alg{A}}_i=\pi(\alg{A}_i)$ under the GNS
  representation of~$\alg{U}$ are free with amalgamation over~$\alg{B}$.
  To this end
  let~$X_j\in\alg{A}$, $1\leq j\leq n$ be an arbitrary finite sequence with
  $\psi_\infty(X_j)=0$
  and let~$h:[n]\to\IN$ be an index function with $h(j)\neq h(j+1)$.
  We have to show that
  $$
  \psi_\infty(X_1^{(h(1))} X_2^{(h(2))}\dotsm X_n^{(h(n))})=0
  $$
  By the assumed faithfulness of~$\phitt$ on~$\alg{B}_0$ it suffices to show that
  $$
  \phitt(\psi_\infty(X_1^{(h(1))} X_2^{(h(2))}\dotsm X_n^{(h(n))})Y)=0
  \qquad\forall Y\in\alg{B}_0
  $$
  and it is enough to consider monomials of the form
  $$
  Y = \psi_\infty(Y_1^{(g(1))} Y_2^{(g(2))}\dotsm Y_m^{(g(m))})
  $$
  with $Y_j\in \alg{A}$ and $g$ an arbitrary index function.
  Indeed, any element of $\alg{B}_0$ is a sum of products
  of elements like this, and for products we have for any~$X\in\alg{B}_0$
  \begin{multline*}
    \phitt(X
  \psi_\infty(Y_1^{(f(1))}Y_2^{(f(2))}\dotsm Y_p^{(f(p))})\,
  \psi_\infty(Z_1^{(g(1))}Z_2^{(g(2))}\dotsm Z_q^{(g(q))}) 
  )\\
  =
    \phitt(X
    \psi_\infty(Y_1^{(f(1))}Y_2^{(f(2))}\dotsm Y_p^{(f(p))}
  Z_1^{(g'(1))}Z_2^{(g'(2))}\dotsm Z_q^{(g'(q))})
  )
  \end{multline*}
  where $g'$ is an index function with $\ker g' = \ker g$ and whose range is
  disjoint from the range of $f$.
  Thus consider
  \begin{align*}
    \phitt(
    &
       \psi_\infty(X_1^{(h(1))}X_2^{(h(2))} \dotsm X_n^{(h(n))})\,
       \psi_\infty(Y_1^{(g(1))} Y_2^{(g(2))}\dotsm Y_m^{(g(m))}))\\
    &= \lim_{N\to\infty}
    \phit(
       \psi_N(X_1^{(h(1))}X_2^{(h(2))} \dotsm X_n^{(h(n))})\,
       \psi_N(Y_1^{(g(1))} Y_2^{(g(2))}\dotsm Y_m^{(g(m))}))\\
    &= \lim_{N\to\infty}
    \phit(
       \psi_N(X_1^{(h(1))}X_2^{(h(2))} \dotsm X_n^{(h(n))}
              Y_1^{(g(1))} Y_2^{(g(2))}\dotsm Y_m^{(g(m))}))\\
\intertext{where we assume without loss of generality that $h$ and $g$ have
    disjoint range,}
    &= \phit(
       X_1^{(h(1))}X_2^{(h(2))} \dotsm X_n^{(h(n))}
       Y_1^{(g(1))} Y_2^{(g(2))}\dotsm Y_m^{(g(m))})\\
    &= \sum_{\substack{\rho_1\leq \ker h \\ \rho_2\leq \ker g}}
        K_{\rho_1\cup\rho_2}^{\exchm}(X_1,\dots,X_n,Y_1,\dots,Y_m)
    ;
  \end{align*}
  by assumption the sum runs over all noncrossing partitions only.
  By Lemma~\ref{lem:GeneralCumulants4:singletonpartition} any noncrossing
  partition $\rho_1\leq\ker h$ contains a singleton, say~$\{j\}$.
  Now
  \begin{align*}
    0 &=\phitt(\psi_\infty(X_j)^*\psi_\infty(X_j)) \\
      &=\lim_{N\to\infty}\phit(\psi_N(X_j)^*\psi_N(X_j)) \\
      &=\lim_{N\to\infty}\phit(\psi_N({X_j^{(1)}}^*X_j^{(2)})) \\
      &=\phit({X_j^{(1)}}^*X_j^{(2)})
  \end{align*}
  and we may apply Proposition~\ref{prop:GeneralCumulant4:WeakSingletonCumulants}
  to every term of the sum to see that it vanishes.
\end{proof}

The main problem is now to prove faithfulness of $\phitt$ on
$\alg{B}_0$.
In the tracial case we may dispose of this problem as follows.
\begin{proof}%
  [End of the proof of Theorem~\ref{thm:GeneralCumulant4:MainTheorem} in the
  tracial case]
  If $\phit$ is a trace, so is~$\phitt$ and its kernel is a two sided ideal.
  Since $\phit$ is faithful, the intersection of $\pi(\alg{U})$ with
  $\ker\phitt$ is trivial
  and therefore $\alg{U}$ is faithfully embedded into the quotient 
  algebra~$\tilde{U}/\ker\phitt$, on which the trace is faithful.
  Now we can apply the arguments of the proof above with $\tilde{U}$
  replaced by the quotient~$\tilde{U}/\ker\phitt$.
\end{proof}
In the non-tracial case there is more work to do, namely
we will show that the extended state $\phitt$ is indeed faithful on
$\alg{B}_0$. To this end we need a very strong law of large numbers for noncrossing
exchangeability systems, which we prove in the next section.

\section{A noncommutative $L^p$-inequality in the case of noncrossing cumulants}
\label{sec:GeneralCumulants4:Lpinequality}

Our aim is to show that $\psi_\infty(X)=0$ if $\phit({X^{(1)}}^*X^{(2)})=0$.
That is, such random variables satisfy a very strong law of large numbers.
We need a combinatorial proof in order to use the combinatorial
information about cumulants, that is, we will use the fact that for a faithful
state $\phit$ we have
\begin{equation}
  \label{eq:GeneralCumulant4:opnorm=limLpnorm}
  \norm{\sum X_i} = \lim_{p\to\infty} \phit( ((\sum X_i)^*(\sum X_i))^p ) )^{1/2p}
\end{equation}
The proof is somewhat in the spirit of~\cite{Pisier:2000:inequality}
where it is shown that the noncommutative $L^p$-norms of so-called
$p$-orthogonal sums (of which our situation is a special case)
can be estimated 
$$
\norm{\sum X_i}_{L^{2p}(\tau)} \leq \frac{3\pi}{2}\,p\,S(X,p)
$$
where
$$
S(X,p) = 
\max
\left\{
  \norm{
    \left(
      \sum X_i^* X_i
    \right)^{1/2}
}
    ,
\norm{
    \left(
      \sum X_i X_i^*
    \right)^{1/2}}
\right\}
$$
However
in order to get something useful out of~\eqref{eq:GeneralCumulant4:opnorm=limLpnorm}
we will need constants which stay bounded as $p$ tends to infinity.
This is related to the question in~\cite[Remark~0.3]{Pisier:2000:inequality} 
whether there are uniform constants for free martingale inequalities,
owing to the fact that the size of the lattice $\NC_n$ of
noncrossing partitions is of order $4^n$, while the size of the lattice of all
partitions $\Pi_n$ is much bigger.
The tracial version in Proposition~\ref{prop:GeneralCumulants4:KhinchinTracial}
gives further evidence for a positive answer to this question.
For our purposes however we need a variant of the inequality
for i.i.d.\ sequences also in the nontracial case.

\begin{Proposition}
  \label{prop:GeneralCumulants4:KhinchinGeneral}
  Assume that a noncrossing exchangeability system $\exch$ satisfies the weak singleton
  condition and has a faithful state.
  Then for any selfajoint random variable~$X$ with~$\phit(X^{(1)}X^{(2)})=0$ the
  interchangeable sequence~$X^{(i)}$ satisfies the inequality
  $$
  \biggnorm{\sum_{i=1}^N X^{(i)}}
  \leq \frac{2\sqrt{N}}{1-\frac{1}{\sqrt{N}}} \norm{X}
  $$
  for every $N\geq2$.
\end{Proposition}
\begin{proof}
  We give three estimates with increasing difficulty and accuracy.
  Roughly the idea is as follows.
  We assume that $X^{(i)}$ are as in the statement of the proposition.
  By faithfulness of~$\phit$, we can use~\eqref{eq:GeneralCumulant4:opnorm=limLpnorm}
  although the ``$L^p$-norm'' associated to $\phit$ is not really a norm.
  We can expand the latter in terms of cumulants:
  \begin{align*}
    \phit(
    \left(
      \sum X^{(i)}
    \right)^p)
    &= \sum_{\pi\in\NC_p}
        K_{\pi}(\sum X^{(i)}) \\
    &= \sum_{\pi\in\NC_p}
       \sum_{\ker h\geq\pi}
        K_{\pi}(X^{(h(1))},X^{(h(2))},\dots,X^{(h(p))}) \\
    &= \sum_{\pi} N^{\abs{\pi}} K_\pi(X)
\intertext{  \emph{First estimate}. Because of the weak singleton condition only partitions without
  singletons contribute. Any such partition has at most $\frac{p}{2}$ blocks
  and therefore the sum is of order $N^{p/2}$ times the number of noncrossing partitions:}
    &\leq N^{p/2} \frac{1}{p+1}\binom{2p}{p} \max_{\pi}\abs{K_\pi(X)}
  \end{align*}
  Each cumulant~$K_\pi$ in turn can be estimated by
  \begin{align*}
    \abs{K_\pi(X)}
    &= \biggabs{\sum_{\sigma\leq\pi} \phi_\sigma(X)\,\mu_{\NC}(\sigma,\pi)}\\
    &\leq \frac{1}{p+1}\binom{2p}{p} \norm{X}^p
    \max_{\sigma,\pi}\abs{\mu_{\NC}(\sigma,\pi)} \\
    &\simeq 16^p \norm{X}^p
  \end{align*}
  Thus by this first rough estimate we obtain the inequality
  $$
  \abs{\phit
    \left(
      (\sum X^{(i)})^p
    \right)}
  \leq 64^p N^{p/2} \norm{X}^p
  $$
  and taking limits
  $$
  \norm{\sum X^{(i)}} \leq 64\sqrt{N} \norm{X}
  .
  $$
  
  \emph{Second estimate}. With a little effort, we can improve on the constant considerably.
  First note that we can evaluate~$a_\pi=\sum_{\sigma\leq\pi}\abs{\mu_{\NC}(\sigma,\pi)}$
  explicitly.
  Since $\mu_{\NC}$ is a multiplicative function, so are $\abs{\mu_{\NC}}$
  and $a=\abs{\mu_{\NC}}\boxstar\zeta$.
  By applying the Kreweras complementation map we have
  \begin{equation}
    \label{eq:GeneralCumulant4:an=sumabsmu}
    a_n
    =\sum_{\sigma\in\NC_n}\abs{\mu_{\NC}(\sigma,\hat1_n)}
    =\sum_{\sigma\in\NC_n}\abs{\mu_{\NC}(\hat0_n,\sigma)}
  \end{equation}
  i.e., $a=\abs{\mu_{\NC}}\boxstar\zeta$ and we can
  use~\eqref{eq:GeneralCumulant4:fboxstarzeta=g}.
  The characteristic series of $\abs{\mu_{\NC}}$ is 
  $$
  \phi_{\abs{\mu}} = \frac{1}{2}(1-\sqrt{1-4z})
  $$
  and $\phi_a(z)$ satisfies the equation
  $$
  \frac{1}{2}
  \left(
    1-\sqrt{1-4z(1+\phi_a(z))}
  \right)
  = \phi_a(z)
  .
  $$
  Together with the condition $\phi_a(0)=0$ this yields the solution
  $$
  \phi_a(z) = \frac{1}{2}(1-z-\sqrt{1-6z+z^2})
  = z + 2z^2 +6z^3 + 22z^4 +\dots
  .
  $$
  This is the generating function of the ``large Schr{\"o}der numbers''
  \cite{Stanley:1999:Enumerative2,Deutsch:2001:bijective};
  they show up in a similar context in~\cite{Dykema:2005:multilinear}.

  We have to estimate
  $$
  \biggabs{
  \sum_{\pi\in\NC_p^{\geq2}} N^{\abs{\pi}} K_\pi(X)
  }
  \leq \sum_{\pi\in\NC_p^{\geq2}} N^{\abs{\pi}} a_\pi \norm{X}^p
  $$
  where $\NC_p^{\geq2}$ is the set of noncrossing partitions without singletons.
  The sequence
  $$
  b_n =  \sum_{\pi\in\NC_n^{\geq2}} N^{\abs{\pi}} a_\pi 
  $$
  is the characteristic sequence of the convolution of the multiplicative 
  function~$N\cdot\bub{a}$ with characteristic sequence~$(N\bub{a}_n)_n$
  with the~$\zeta$-function, where
  $$
  \bub{a}_n=
  \begin{cases}
    0 & n=1\\
    a_n& n\geq 2
  \end{cases}
  $$
  and
  $$
  \phi_{\bub{a}}(z)
  = \frac{1}{2}
    \left(
      1-3z-\sqrt{1-6z+z^2}
    \right)
  $$
  The characteristic series~$\phi_b(z)$ can be found by yet another appeal
  to~\eqref{eq:GeneralCumulant4:fboxstarzeta=g}, namely it satisfies the equation
  $$
  N\phi_{\bub{a}}(z(1+\phi_b(z)))=\phi_b(z)
  $$
  and the relevant solution is
  $$
  \phi_b(z) = 
  \frac{2(N+1)}{N+2+3Nz+N\sqrt{1-6z+(1-8N)z^2}} - 1
  .
  $$
  The dominant singularity comes from the radical $1-6z+(1-8N)z^2$.
  The zeros of the latter are $\frac{1}{8N-1}
  \left(
    \pm\sqrt{2(N+1)}-3
  \right)
  $
  and therefore
  $$
  b_n \sim
  \left(
    \frac{8N-1}{2\sqrt{2(N+1)}-3}
  \right)^n
  ;
  $$
  it follows that
  $$
  \norm{\sum X^{(i)}} \leq 
  \frac{8N-1}{2\sqrt{2(N+1)}-3} \norm{X}
  .
  $$
  The constant tends to~$2\sqrt{2}$ as~$N\to\infty$, which is not bad, as the
  best possible constant is~$2$.
  
  \emph{Third estimate}.
  With even some more effort, one can obtain the optimal constant (at least as
  $N\to\infty$)   as follows.
  The previous estimate was done using the numbers $a_n$
  from~\eqref{eq:GeneralCumulant4:an=sumabsmu} 
  and we neglected the fact that for the calculation of the cumulants $K_\pi(X)$
  partitions with singletons do not contribute.
  Thus it will be more accurate to work with the numbers
  \begin{equation}
    \label{eq:GeneralCumulant4:an=sumabsmuwosingletons} 
    \tilde{a}_n = \sum_{\pi\in\NC_n^{\geq2}} \abs{\mu_{\NC}(\pi,\hat1_n)}
  \end{equation}
  which constitute the characteristic sequence of the multiplicative function
  $\bub{\zeta}\boxstar\abs{\mu_{\NC}}$ where
  $$
  \bub{\zeta}_n=
  \begin{cases}
    0 & n=1\\
    1 & n\geq 2
  \end{cases}
  $$
  is the Zeta function on the poset of noncrossing partitions without singletons.
  This convolution can be carried out with the aid
  of~\eqref{eq:GeneralCumulant4:Ffboxstarg=FfFg}. The ``Fourier transforms''
  of the functions
  $$
  \phi_{\bub\zeta}(z) = \sum_{n=2}^\infty z^n = \frac{z^2}{1-z}
  \qquad\text{and}\qquad
  \phi_{\abs{\mu}}(z) = \frac{1}{2}
  \left(
    1-\sqrt{1-4z}
  \right)
  $$
  are
  $$
  \alg{F}_{\bub{\zeta}}
  =\frac{\pm\sqrt{z^2+4z}-z}{2z}
  \qquad\text{and}\qquad
  \alg{F}_{\abs{\mu}}(z) = 1-z
  $$
  respectively.
  Therefore
  $$
  \alg{F}_{\bub{\zeta}\boxstar\abs{\mu}}(z) =
  \frac{\pm\sqrt{z^2+4z}-z}{2z}
  \,
  (1-z)
  $$
  i.e., $y=y(z)=\phi_{\bub{\zeta}\boxstar\abs{\mu}}(z)$ satisfies the
  algebraic equation
  $$
  y(1-y)(1-y-z)=z^2
  .
  $$
  We are interested in the asymptotics of the numbers
  $$
  \tilde{b}_n = \sum_{\pi\in\NC_n} N^{\abs{\pi}} \tilde{a}_\pi
  $$
  whose generating function can be determined
  by~\eqref{eq:GeneralCumulant4:fboxstarzeta=g}, namely
  $$
  N\phi_{\tilde{a}}(z(1+\phi_{\tilde{b}}(z))) = \phi_{\tilde{b}}(z)
  $$
  Thus $x=x(z)=\phi_{\tilde{b}}(z)$ satisfies the equations
  $$
  N \phi_{\tilde{a}}(z(1+x)) = x
  \qquad
  \frac{x}{N}(1-\frac{x}{N})(1-\frac{x}{N}-z(1+x)) = z^2(1+x)^2
  $$
  therefore $x=x(z)$ is the solution of the equation
  \begin{equation}
    \label{eq:GeneralCumulant4:equationforphib}
  g(x,z)=
  \frac{x}{N}\, \left(1 - \frac{x}{N}\right)\, \left(1 - \frac{x}{N} - z\,
  \left(x + 1\right)\right) - z^2\, \left(x + 1\right)^2 
  =0
  \end{equation}
  If $x(z)$ has a singularity at~$z$, then both $g(x,z)=0$ and $\partial_x
  g(x,z)=0$
  and therefore~$z$ is a zero of the
  resultant 
  \begin{multline*}
    Res(g(x,z),\partial_zg(x,z)) =
    \frac{1}{N^9}   (N\, z + 1)\, (N + 1)^2\, z^2 
\\     \times
    \left(8\, z - 2\, N\, z + 32\, z^2 - 4\, z^3 + 26\, N\, z^2 + 16\, N\,
      z^3 - N^2\, z^2 + 10\, N^2\, z^3 - N^2\, z^4 + 4\, N^3\, z^4 - 5
    \right)
    ,
  \end{multline*}
  cf.~\cite{CoxLittleOShea:1998:using,FlajoletSedgewick:2001:chapter8}.
  By Pringsheim's theorem we know that the dominant singularity is positive and
  therefore it must be a root of the last factor
  $$
  r(z) = -5 + (8-2N)\,z + (32+26N-N^2)\,z^2 + (-4+16N+10N^2)\,z^3
  +(-N^2+4N^3)\,z^4
  $$
  We claim that $r(z)\neq0$ for $0\leq z\leq\frac{1}{2\sqrt{N}}
  \left(
    1-\frac{1}{\sqrt{N}}
  \right)
  .
  $
  Indeed, let~$z=\frac{\alpha}{2\sqrt{N}}
  \left(
    1-\frac{1}{\sqrt{N}}
  \right)$ with
  $0\leq\alpha\leq1$,
  then it is tedious but not difficult to verify that
  \begin{align*}
    r(\frac{\alpha}{2\sqrt{N}}&
  \left(
    1-\frac{1}{\sqrt{N}}
  \right))\\
  &=N
  \left(
    \frac{\alpha^4}{4}-\frac{\alpha^2}{4}
  \right)
  +
  N^{1/2}
  \left(
    -\alpha+\frac{\alpha^2}{2}+\frac{5}{4}\,\alpha^3-\alpha^4
  \right)\\
  &\phantom=
  - 5 + \alpha+\frac{25}{4}\,\alpha^2-\frac{15}{4}\,\alpha^3+\frac{23}{16}\,\alpha^4
  + N^{-1/2} (4\alpha-13\alpha^2+\frac{23}{4}\,\alpha^3-\frac{3}{4}\,\alpha^4)\\
  &\phantom=
  + N^{-1} (-4\alpha+\frac{29}{2}\,\alpha^2-\frac{29}{4}\,\alpha^3-\frac{1}{8}\,\alpha^4)
  + N^{-3/2} (-16\alpha^2+\frac{11}{2}\,\alpha^3+\frac{1}{8}\,\alpha^4)\\
  &\phantom=
  + N^{-2} (8\,\alpha^2-\frac{1}{2}\,\alpha^3-\frac{1}{16}\,\alpha^4)
  - N^{-5/2} \frac{3}{2}\,\alpha^3
  + N^{-3} \frac{1}{2}\,\alpha^3\\
  &=
  -\frac{N}{4}\,\alpha^2(1-\alpha^2)
  -
  N^{1/2}
  \left(
    \frac{1}{4}(\frac{1}{4}-(\alpha-\frac{1}{2})^2)+(1-\alpha) + \alpha(\alpha-\frac{1}{2})^2
  \right)\\
  &\phantom=
  -
  (1 - N^{-1/2})
  \left(
    \frac{1}{16} + (1-\alpha)
    \left(
      \frac{47}{16} + \frac{29}{8}\,\alpha + 2(1-\alpha)^2+\frac{23}{16}(1-(1-\alpha)^3)
    \right)
  \right)\\
  &\phantom=
  -
  (N^{-1/2}-N^{-1})
  \biggl(
    \frac{877}{256} + \frac{3}{32}\,(1-\alpha) 
    +\frac{27}{8}
    \left(
      \frac{1}{8}-
      \left(
        \alpha-\frac{1}{2}
      \right)^3
    \right)
\\ &\phantom{=====================}
    +
    \left(
      \alpha-\frac{1}{2}
    \right)^2
    \left(
      \frac{685}{256}
      +\frac{11}{16}\,
      \left(
        \frac{1}{4}
        -
        \left(
          \alpha-\frac{1}{2}
        \right)^2
      \right)
    \right)
  \biggr)\\
  &\phantom=  
  -N^{-1}
  \left(
    \frac{15}{16}
    + 3(1-\alpha)(1-(1-\alpha)^2)
    + \frac{37}{8} (1-\alpha)^2
    -\frac{9}{16} (1-\alpha)^4
  \right)\\
  &\phantom=  
  -(N^{-3/2}-N^{-2})
  \alpha^2(16-\frac{11}{2}\,\alpha-\frac{1}{4}\,\alpha^2)
  -N^{-2} \alpha^2(8-5\alpha-\frac{3}{16}\,\alpha^2)\\
  &\phantom=  
  -N^{-5/2} (3-N^{-1/2})\frac{1}{2}\,\alpha^3
  \end{align*}
  which is strictly negative for $0\leq\alpha\leq1$.
  Therefore asymptotically as~$n$ tends to infinity we have
  $$
  \tilde{b}_n\leq
  \left(
    \frac{2\sqrt{N}}{1-\frac{1}{\sqrt{N}}} 
  \right)^n
  .
  $$
\end{proof}
Using Remark~\ref{rem:GeneralCumulant4:iid=interchangeable} we obtain the following Corollary.
\begin{Corollary}
  \label{cor:GeneralCumulants4:KhinchinGeneral}
  Assume that a noncrossing exchangeability system
  $\exch=(\alg{U},\phit,\alg{J})$ satisfies (\eWSC{}) and has a faithful state.
  Let $I\subseteq \IN$ be a finite index set and $(I_j)$ a sequence of disjoint
  index sets of the same cardinality as $I$, cf.\  
  Remark~\ref{rem:GeneralCumulant4:iid=interchangeable}.
  Let~$X\in\alg{A}_I$ be a selfadjoint polynomial
  with~$\phit(X^{(I_1)}X^{(I_2)})=0$,
  then the
  interchangeable sequence~$X^{(I_j)}$ satisfies the inequality
  $$
  \biggnorm{\sum_{j=1}^N X^{(I_j)}}
  \leq \frac{2\sqrt{N}}{1-\frac{1}{\sqrt{N}}} \norm{X}
  $$
  for every $N\geq2$.
  In particular, $\psi_N(X^{(I)})$ is $\Ord(1/\sqrt{N})$ and converges to zero
  as $N$ tends to infinity.
\end{Corollary}
\begin{Remark}
  This is the only place where (\eWSC{}) is needed rather than (\WSC{}).
  While (\eWSC{}) holds in all examples known to us, we were not able to decide
  whether it follows from (\WSC{}).
\end{Remark}

\begin{proof}[End of the proof of Theorem~\ref{thm:GeneralCumulant4:MainTheorem} in the
  nontracial case]
  It remains to prove faithfulness of the state~$\phit$ on~$\alg{B}_0$.
  Let~$X\in\alg{B}_0$ such that $\phit(X^*X)=0$,
  i.e.~$X=\psi_\infty(W)$, where~$W\in\alg{U}_0$ is some polynomial, say $W\in
  A_I$ for some finite index set $I$, and let $I_1$ and $I_2$ be disjoint
  copies of $I$, cf.\ Remark~\ref{rem:GeneralCumulant4:iid=interchangeable}.
  By assumption
  \begin{align*}
    0&= \phit(\psi_\infty(W)^*\psi_\infty(W))\\
     &= \phit(\psi_\infty({W^{(I_1)}}^* W^{(I_2)}))\\
     &= \phit({W^{(I_1)}}^* W^{(I_2)})\\
  \end{align*}
  and by Corollary~\ref{cor:GeneralCumulants4:KhinchinGeneral} this implies
  that~$\psi_\infty(W)=\lim_{N\to\infty}\psi_N(W)=0$.
\end{proof}

\section{Weak freeness}
\label{sec:GeneralCumulant4:WeakFreeness}

In this section we discuss the notion of weak freeness, which together with 
(\eWSC{}) implies vanishing of crossing cumulants.
\begin{Definition}
  \label{def:GeneralCumulant4:WeakFreeness}
  Let~$\exchm=(\alg{U},\phit,\alg{J})$ be an exchangeability system for
  a noncommutative probability   space~$(\alg{A},\phi)$. 
  For an index set $I\subseteq\IN$ denote $\alg{A}_I$ the subalgebra of
  $\alg{U}$ generated by $(\alg{A}_i)_{i\in I}$.
  We say that~$\exchm$ satisfies \emph{weak freeness}
  if
  $$
  \phit(X_1X_2\dotsm X_n)=0
  $$
  whenever~$\phit({X_j^{(1)}}^*X_j^{(2)})=0$, $I_k$ are disjoint index sets and
  $X_j\in\alg{A}_{(I_{i_j})}$ with~$i_j\neq i_{j+1}$ for every~$j$.
  Here $X_j^{(1)}$ and $X_j^{(2)}$ refer to copies of $X_j$ in $A_{I'_{i_j}}$
  and $A_{I''_{i_j}}$, where $I'_{i_j}$ and $I''_{i_j}$ are disjoint copies
  of $I_{i_j}$.
\end{Definition}

It will be convenient to adapt the exchangeability system as indicated in
Remark~\ref{rem:GeneralCumulant4:iid=interchangeable}.
Decompose~$\IN$ into an infinite union of disjoint copies of itself
$\IN=\bigcup_{j=0}^\infty I_j$. Then relabel the indices and consider
the exchangeability system with embeddings
$\iota_{ij}:\alg{A}\to\alg{A}_{i,j}\subseteq\alg{U}$, $i,j\in\IN$.
Thus~$\alg{U}$ is also an exchangeability system for
$\tilde{A}=\bigvee_{j\in\IN}\alg{A}_{0j}$
and we will work with this interpretation in this section, i.e.,
our random variables~$X$ are elements of~$\tilde{\alg{A}}$ and
$X^{(i)}$ are elements of $\tilde{\alg{A}}_i=\bigvee_{j\in\IN}\alg{A}_{ij}$.
Thus if $X=X_1^{(0,j_1)} X_2^{(0,j_2)}\dotsm X_n^{(0,j_n)}$,
then $X^{(i)}=X_1^{(i,j_1)} X_2^{(i,j_2)}\dotsm X_n^{(i,j_n)}$.
The weak freeness condition of
Definition~\ref{def:GeneralCumulant4:WeakFreeness} can be rephrased more 
clearly as follows, namely
$$
\phit(X_1^{(i_1)}X_2^{(i_2)}\dotsm X_n^{(i_n)})=0
$$
whenever $X_j\in\tilde{\alg{A}}$ with $\phit({X_j^{(1)}}^*X_j^{(2)})=0$ and
$i_j\neq i_{j+1}$ for all $1\leq j\leq n-1$.

We need this regrouping in order to define an asymptotic conditional
expectation which is used to transfer proofs from the amalgamated free
situation. As in Section~\ref{subsec:GeneralCumulant4:ConditionalExpectations},
we define symmetrizing maps
\begin{align*}
  \psi_N:\tilde{\alg{A}} & \to\tilde{\alg{A}}\\
      X &\mapsto \frac{1}{N!} \sum_{\sigma\in\SG_N} \sigma(X)
\end{align*}
which will allow us to construct asymptotically $\psi$-centered random variables.
\begin{Lemma}
  \label{lem:GeneralCumulant4:AsymptoticMixedCumulantsVanish}
  Let $X_j\in\tilde{\alg{A}}$ be polynomials, that is, linear combinations of
  elements of the form
  $$
  Z_1^{(0,j_1)} Z_2^{(0,j_2)}\dotsm Z_m^{(0,j_m)}
  $$
  with $Z_j\in\alg{A}$.
  Then
  $$
  K_\pi(X_1,\dots,X_{k-1},\psi_N(X_k),X_{k+1},\dots,X_n)\xrightarrow[N\to\infty]{} 0
  $$
  unless $\{k\}$ is a singleton of~$\pi$.
  In other words, $\psi_N(X_k)$ is asymptotically independent from the rest.
\end{Lemma}
\begin{proof}
  Indeed
  $$
  K_\pi(X_1,\dots,X_{k-1},\psi_N(X_k),X_{k+1},\dots,X_n)
  =\frac{1}{N!} \sum_{\sigma\in\SG_N} K_\pi(X_1,\dots,X_{k-1},\sigma(X_k),X_{k+1},\dots,X_n)
  $$
  Let $s$ be the maximal superscript appearing in the polynomials~$X_j$.
  If $\sigma$ maps each $1\leq j\leq s$ to some index strictly greater than
  $s$, then $\sigma(X_k)$ is independent from the other~$X_j$ and by assumption
  the cumulant vanishes. The number of permutations~$\sigma$ of this type is
  $(N-s)(N-s-1)\dotsm(N-2s)\cdot(N-s)!$, i.e., almost all permutations, because
  the ratio of the rest is
  $$
  \frac{N!-(N-s)(N-s-1)\dotsm(N-2s)\cdot(N-s)!}{N!}\xrightarrow[N\to\infty]{}0.
  $$
\end{proof}

When calculating cumulants of elements $X_j\in \tilde{\alg{A}}$,
we can thus replace $X_j$ by $X_j-\psi_N(X_j)$ 
for each non-singleton index $j$ and then let $N$ tend to infinity.
Replacing $X_j$ by $X_j-\psi_N(X_j)$ allows to apply the weak singleton
condition,
because due to the permutation invariance of the state we have
$$
\phit({\sigma(X)^{(1)}}^*X^{(2)}) = \phit({X^{(1)}}^*X^{(2)}) 
$$
and therefore
\begin{align*}
  \phit((X^{(1)}-\psi_N(X)^{(1)})^* & (X^{(2)}-\psi_N(X)^{(2)})) \\
  &= \phit(X^{(1)}{}^*(X)^{(2)})
     -
     \frac{1}{N!}\sum_\sigma \phit(\sigma(X)^{(1)}{}^*X^{(2)}) \\
  &\phantom{==}  -
     \frac{1}{N!}\sum_\sigma \phit(X^{(1)}{}^*\sigma(X)^{(2)})
     +
     \frac{1}{(N!)^2}\sum_{\sigma,\tau} \phit(\sigma(X)^{(1)}{}^*\tau(X)^{(2)})\\
  &= 0
\end{align*}

\begin{Theorem}
  Let $\exchm=(\alg{U},\phit,\alg{J})$
  be an exchangeability system for a noncommutative probability space
  $(\alg{A},\phi)$ with faithful state~$\phit$ such that 
  both weak freeness and the weak singleton condition holds.
  Then crossing cumulants vanish.
  In particular, $\exchm$ can be embedded into an amalgamated free product.
\end{Theorem}

\begin{proof}
  The proof consists of three parts by reducing an arbitrary crossing partition 
  to an alternating partition  
  without singletons.

  We will use the following terminology.
  A partition~$\pi\in\Pi_n$ is called \emph{alternating} 
  if $i\not\sim_\pi i+1$ for all $1\leq i\leq n-1$, that is, adjacent elements are
  in different blocks of $\pi$. For non-alternating partitions we denote by
  $$
  \cn(\pi)=\#\{(k,k+1) : k\sim_\pi k+1,1\leq k\leq n-1\}
  $$
  the \emph{number of connected neighbours} of~$\pi$. Clearly~$\pi$ is alternating
  if and only if~$\cn(\pi)=0$.
  
  \emph{Step 1. Alternating partitions without singletons.}
  First assume that $\pi$ is an alternating partition without singletons.
  Let~$\eps>0$.
  Then by Lemma~\ref{lem:GeneralCumulant4:AsymptoticMixedCumulantsVanish}
  we may find $N>0$ such that
  $$
  K_\pi(X_1,\dots,X_n) =  K_\pi(X_1-\psi_N(X_1),\dots,X_n-\psi_N(X_n)) +R_n
  $$
  with error term  $\abs{R_n}<\eps$.
  Now a look at the moment-cumulant formula
  $$
  K_\pi(X_1-\psi_N(X_1),\dots,X_n-\psi_N(X_n))
  = \sum_{\sigma\leq\pi}
     \phi_\sigma(X_1-\psi_N(X_1),\dots,X_n-\psi_N(X_n))
     \,\mu(\sigma,\pi)
  $$
  shows that the sum runs over alternating partitions and by weak freeness every
  term vanishes.

  \emph{Step 2. Reducing everything to alternating partitions.}
  The aim is now to express an arbitrary cumulant in terms of alternating ones.
  We will use the product formula of Leonov and Shiryaev from
  Proposition~\ref{prop:GeneralCumulants:Productformula} to reduce the number
  of connected neighbours.
  Consider a partition~$\pi\in\Pi_n$ with crossings and $\cn(\pi)>0$.
  Pick an arbitrary element $k$ with $k\sim_\pi k+1$. We will express $K_\pi$
  as a sum of cumulants $K_\rho$ s.t.\  $k\not\sim_\rho k+1$ and
  $\cn(\rho)<\cn(\pi)$.
  Let $\hat\pi=\pi/[k=k+1]\in\Pi_{n-1}$ be the partition obtained from~$\pi$ by
  identifying~$k$ and~$k+1$. Then by
  Proposition~\ref{prop:GeneralCumulants:Productformula} we have
  with~$\nu=\{\{1\},\{2\},\dots,\{k,k+1\},\dots,\{n\}=
\begin{picture}(110,8.4)(1,0)
  \put(10,0){\line(0,1){8.4}}
  \put(20,0){\line(0,1){8.4}}
  \put(30,0){\ldots}
  \put(50,0){\line(0,1){8.4}}
  \put(60,0){\line(0,1){8.4}}
  \put(70,0){\line(0,1){8.4}}
  \put(80,0){\line(0,1){8.4}}
  \put(90,0){\ldots}
  \put(110,0){\line(0,1){8.4}}
  \put(10,8.4){\line(1,0){0}}
  \put(20,8.4){\line(1,0){0}}
  \put(50,8.4){\line(1,0){0}}
  \put(60,8.4){\line(1,0){10}}
  \put(80,8.4){\line(1,0){0}}
  \put(110,8.4){\line(1,0){0}}
\end{picture}
$ the decomposition
\begin{align*}
  K_{\hat\pi}(X_1,\dots,X_k X_{k+1},\dots, X_n)
  &= \sum_{\rho\lor \nu=\pi}   K_{\rho}(X_1,\dots,X_k, X_{k+1},\dots, X_n)\\
  &= K_{\pi}(X_1,\dots, X_n) 
     +\sum_{\substack{\rho\lor \nu=\pi\\ \rho<\pi}} K_{\rho}(X_1,\dots, X_n) 
\end{align*}
  Each contributing partition~$\rho$ in the second sum is obtained from~$\pi$ by splitting
  the block containing $k$ and $k+1$ into two in such a way that $k$ and $k+1$ are separated.
  In particular, $\cn(\rho)\leq\cn(\pi)-1$ and $\cn(\hat\pi)=\cn(\pi)-1<\cn(\pi)$ as well.
  To conclude, we have
  $$
  K_{\pi}(X_1,\dots, X_n) 
  =  K_{\hat\pi}(X_1,\dots,X_k X_{k+1},\dots, X_n)
     -\sum_{\substack{\rho\lor \nu=\pi\\ \rho<\pi}} K_{\rho}(X_1,\dots, X_n) 
  $$  
  and all partitions appearing on the right hand side have less connected neighbours
  than~$\pi$. Repeating this operation finitely many times we end up with a linear combination
  of alternating partitions (possibly involving singletons).

  \emph{Step 3. Getting rid of singletons.}  
  Assume that after step~2 we have arrived at an alternating partition~$\pi$.
  Then for~$N$ large enough, we have by 
  Lemma~\ref{lem:GeneralCumulant4:AsymptoticMixedCumulantsVanish} again
  $$
  K_\pi(X_1,\dots,X_n)\approx  K_\pi(\tilde{X}_1,\dots,\tilde{X}_n)
  $$
  where
  $$
  \tilde{X}_j = 
  \begin{cases}
    X_j               & \text{if $\{j\}$ is a singleton of $\pi$}\\
    X_j - \psi_N(X_j) & \text{if $\{j\}$ is not a singleton of $\pi$}
  \end{cases}
  $$
  We may therefore assume without loss of generality that
  $\psi({X_j^{(1)}}^*X_j^{(2)})=0$ for the non-singleton indices~$j$.
  If the singleton entries~$k$ satisfy $\psi({X_k^{(1)}}^*X_k^{(2)})=0$ as well,
  then we may proceed as in step one.
  If however  there are singletons for which this is not the case,
  we may eliminate them as follows.
  Let~$\{k\}$ be the first of these critical singletons, then
  we may write the cumulant as
  $$
  K_\pi(X_1,\dots,X_n)
  = K_\pi(X_1,\dots,X_k-\psi_N(X_k),\dots,X_n)
  + K_\pi(X_1,\dots,\psi_N(X_k),\dots,X_n)
  $$
  The first term has one critical singleton less than the left hand side
  and the second term can be treated with the product formula as follows.
  Let again $\hat\pi = \pi/[k=k+1]\in\Pi_{n-1}$ be the partition obtained
  from~$\pi$ by identifying~$k$ with~$k+1$. 
  Moreover put~$\nu=\{\{1\},\{2\},\dots,\{k,k+1\},\dots,\{n\}=
  \begin{picture}(110,8.4)(1,0)
    \put(10,0){\line(0,1){8.4}}
    \put(20,0){\line(0,1){8.4}}
    \put(30,0){\ldots}
    \put(50,0){\line(0,1){8.4}}
    \put(60,0){\line(0,1){8.4}}
    \put(70,0){\line(0,1){8.4}}
    \put(80,0){\line(0,1){8.4}}
    \put(90,0){\ldots}
    \put(110,0){\line(0,1){8.4}}
    \put(10,8.4){\line(1,0){0}}
    \put(20,8.4){\line(1,0){0}}
    \put(50,8.4){\line(1,0){0}}
    \put(60,8.4){\line(1,0){10}}
    \put(80,8.4){\line(1,0){0}}
    \put(110,8.4){\line(1,0){0}}
  \end{picture}
  $
  and let $\tilde{\pi}=\pi\lor\nu$ be the partition obtained from $\pi$ by
  adjoining the singleton $\{k\}$ to the block containing $k+1$.
  Then we have by Proposition~\ref{prop:GeneralCumulants:Productformula} again
  \begin{multline*}
  K_{\hat\pi}(X_1,X_2,\dots,\psi_N(X_k) X_{k+1},\dots,X_n) 
  = K_{\tilde\pi}(X_1,X_2,\dots,\psi_N(X_k), X_{k+1},\dots,X_n) \\
  + \sum_{\substack{\rho\lor\nu=\tilde\pi\\ \rho<\tilde\pi}} 
     K_{\rho}(X_1,X_2,\dots,\psi_N(X_k), X_{k+1},\dots,X_n)
  \end{multline*}
  $K_{\tilde\pi}$ vanishes asymptotically
  by Lemma~\ref{lem:GeneralCumulant4:AsymptoticMixedCumulantsVanish}
  and so do all $K_\rho$ in which $k\sim_\pi k+1$, and hence
  the only nontrivial term on the right hand side is the cumulant indexed
  by $\rho=\pi$, because this is the only one in which $k$ is a singleton.
  Thus
  $$
  K_{\pi}(X_1,X_2,\dots,\psi_N(X_k), X_{k+1},\dots,X_n) 
  \approx
  K_{\hat\pi}(X_1,X_2,\dots,\psi_N(X_k) X_{k+1},\dots,X_n) 
  $$
  and~$\hat\pi$ has one singleton less than~$\pi$.
  Repeating this procedure we end up with a linear combination of
  alternating partitions without singletons and step one of the proof applies.
\end{proof}

\section{Appendix: Free $L^p$ Inequalities}

We consider now the tracial version of
Proposition~\ref{prop:GeneralCumulants4:KhinchinGeneral}.
In this case H{\"o}lder's inequality is available and we can get estimates
in terms of the $L^p$-norms.
\begin{Proposition}
  \label{prop:GeneralCumulants4:KhinchinTracial}
  Assume that a noncrossing exchangeability
  system~$\exch=(\alg{U},\tau,\alg{J})$ is tracial and faithful.
  Then for any sequence of $\exch$-independent random variables~$X_i$
  with~$\tau({X_i^{(1)}}^*X_i^{(2)})=0$ the inequality
  $$
  \biggnorm{\sum_{i=1}^N X_i}_{L^{2p}(\tau)}
  \leq C_{2p} S(X,2p)
  $$
  holds with~$C_{2p}\leq \frac{3\pi}{4z_0}\simeq 9.85859$ as~$p\to\infty$
  where~$z_0$ is computed in~\eqref{eq:GeneralCumulant4:z0pisierconstant}.
\end{Proposition}
\begin{proof}
  We expand the $L^p$-norm in terms of cumulants:
  \begin{align*}
  \tau
  \left(
    \left(
      \sum X_i^* X_j
    \right)^p
  \right)
  &= \sum_{\pi\in\NC_{2p}} 
      K_\pi(\sum X_i^*,\sum X_i,\dots,\sum X_i^*,\sum X_i)\\
  \end{align*}
  and for each~$\pi$, we estimate the cumulant $K_\pi$.
  Because of the weak singleton condition, only partitions without singletons
  are involved in the sums.
  \begin{align*}
    K_\pi&(\sum X_i^*,\sum X_i,\dots,\sum X_i^*,\sum X_i)\\
    &= \sum_{\ker h\geq\pi}
        K_\pi( X_{h(1)}^*, X_{h(2)},\dots, X_{h(2p-1)}^*, X_{h(2p)})\\
    &= \sum_{\ker h\geq\pi}
        \sum_{\sigma\leq\pi} 
         \tau_\sigma( X_{h(1)}^*, X_{h(2)},\dots, X_{h(2p-1)}^*, X_{h(2p)})
         \,\mu_{\NC}(\sigma,\pi)
  \end{align*}
  Now for fixed~$\sigma$
  we have~\cite[Sublemma~3.3 and Lemma~3.4]{Pisier:2000:inequality}
  $$
  \abs{\sum_{\ker h\geq\pi} \tau_\sigma(X_{h(1)}X_{h(2)}\dotsm X_{h(2p)}) }
  \leq \norm{\sum \lambda(g_i)\ox X_i}_{2p}^{2p}
  \leq \left(
         \frac{3\pi}{4}  
       \right)^{2p} 
       S(X,2p)
  $$
  where~$\lambda(g_i)$ is the left regular representation of the generators of
  the free group.
  Indeed, we can find a suitable discrete group~$G$ and elements
  $F_1$,\ldots,$F_{2p}$ in $L^p(\tau_G\ox \tau)$ such that 
  $$
  \norm{F_k}_{2p} = \norm{\sum \lambda(g_i)\ox X_i}_{2p}
  $$
  and
  $$
  \sum_{\ker h\geq\pi} \tau_\sigma(X_{h(1)}X_{h(2)}\dotsm X_{h(2p)}) 
  = \tau_G\ox\tau(F_1F_2\dotsm F_{2p})
  .
  $$
  The construction is done as in the proof
  of~\cite[Sublemma~3.3]{Pisier:2000:inequality},
  with the slight modification that $X_i$ is replaced by $X_i^{(\sigma(k))}$ when
  constructing~$F_k$.
  With these preparations we continue similarly as before
  $$
  \tau
  \left(
    \left(
      \sum X_i^* X_j
    \right)^p
  \right)
  \leq \sum_{\pi\in\NC^{\geq2}}
        \sum_{\substack{\sigma\in\NC_{2p}^{\geq2}\\\sigma\leq\pi}}
         \abs{\mu_{NC}(\sigma,\pi)}
         \left(
           \frac{3\pi}{4}  
         \right)^{2p} 
         S(X,2p)         
  $$
  We know the asymptotics of
  $$
  \tilde{b}_n = \sum_{\pi\in\NC^{\geq2}}
        \sum_{\substack{\sigma\in\NC_{2p}^{\geq2}\\\sigma\leq\pi}}
         \abs{\mu_{NC}(\sigma,\pi)}
  $$
  from its generating function, which satisfies 
  equation~\eqref{eq:GeneralCumulant4:equationforphib} with~$N=1$.
  Its dominant singularity is a zero of the resultant
  $$
  r(z) = 3z^4+22z^3+57z^2+6z-5
  $$
  and the singularity in question is
  \begin{equation}
    \label{eq:GeneralCumulant4:z0pisierconstant}
    z_0=
    -\frac{11}{6} + \frac{1}{6}{\sqrt{7 + \gamma}}  + {\sqrt{14 - \gamma 
          +
          \frac{992}
          {{\sqrt{7 + \gamma}}}}}
     \simeq 0.238999
   \end{equation}
   where
   $$
   \gamma=        9\,{\left( 207 - 48\,{\sqrt{3}} \right) }^{\frac{1}{3}} +
        9\,{\left( 207 + 48\,{\sqrt{3}} \right) }^{\frac{1}{3}} 
   ;
   $$
   consequently
   $
   \tilde{b}_n \simeq z_0^{-n}
   $
   as $n\to\infty$.
\end{proof}

\subsection*{Acknowledgements}
We are grateful to several anonymous referees who read this paper
during its odyssey for pointing out gaps in proofs
and other improvements.


\begin{thebibliography}{CLO98}

\bibitem[AL93]{AccardiLu:1993:continuous}
L.~Accardi and Y.~G. Lu, \emph{A continuous version of de {F}inetti's theorem},
  Ann. Probab. \textbf{21} (1993), 1478--1493. \MR{94m:60016}

\bibitem[BLS96]{BozejkoLeinertSpeicher:1996:convolution}
Marek Bo{\.z}ejko, Michael Leinert, and Roland Speicher, \emph{Convolution and
  limit theorems for conditionally free random variables}, Pacific J. Math.
  \textbf{175} (1996), 357--388.

\bibitem[BS96]{BozejkoSpeicher:1996:interpolations}
Marek Bo{\.z}ejko and Roland Speicher, \emph{Interpolations between bosonic and
  fermionic relations given by generalized {B}rownian motions}, Math. Z.
  \textbf{222} (1996), 135--159.

\bibitem[CLO98]{CoxLittleOShea:1998:using}
David Cox, John Little, and Donal O'Shea, \emph{Using algebraic geometry},
  Graduate Texts in Mathematics, vol. 185, Springer-Verlag, New York, 1998.
  \MR{99h:13033}

\bibitem[CT78]{ChowTeicher:1978:probability}
Yuan~Shih Chow and Henry Teicher, \emph{Probability theory}, Springer-Verlag,
  New York, 1978. \MR{80a:60004}

\bibitem[Deu01]{Deutsch:2001:bijective}
Emeric Deutsch, \emph{A bijective proof of the equation linking the
  {S}chr{\"o}der numbers, large and small}, Discrete Math. \textbf{241} (2001),
  235--240. \MR{2002h:05011}

\bibitem[DF80]{DiaconisFreedman:1980:finite}
P.~Diaconis and D.~Freedman, \emph{Finite exchangeable sequences}, Ann. Probab.
  \textbf{8} (1980), no.~4, 745--764. \MR{81m:60032}

\bibitem[Dyk05]{Dykema:2005:multilinear}
Ken Dykema, \emph{{Multilinear function series and transforms in free
  probability theory}}, 2005, Preprint, arXiv:math.OA/0504361.

\bibitem[FS01]{FlajoletSedgewick:2001:chapter8}
P.~Flajolet and R.~Sedgewick, \emph{Analytic combinatorics: functional
  equations, rational, and algebraic functions}, preprint,
  \url{http://pauillac.inria.fr/algo/flajolet/Publications/FlSe01.pdf}, 2001,
  Chapter~8 of forthcoming book ``Analytic Combinatorics''.

\bibitem[HM76]{HudsonMoody:1975:locally}
R.~L. Hudson and G.~R. Moody, \emph{Locally normal symmetric states and an
  analogue of de {F}inetti's theorem}, Z. Wahrscheinlichkeitstheorie und Verw.
  Gebiete \textbf{33} (1975/76), 343--351. \MR{53 \#1280}

\bibitem[Hud81]{Hudson:1981:analogs}
R.~L. Hudson, \emph{Analogs of de {F}inetti's theorem and interpretative
  problems of quantum mechanics}, Found. Phys. \textbf{11} (1981), 805--808.
  \MR{83f:81010}

\bibitem[Kin78]{Kingman:1978:uses}
J.~F.~C. Kingman, \emph{Uses of exchangeability}, Ann. Probability \textbf{6}
  (1978), 183--197. \MR{58 \#13238}

\bibitem[Kre72]{Kreweras:1972:partitions}
G.~Kreweras, \emph{Sur les partitions non crois{\'e}es d'un cycle}, Discrete
  Math. \textbf{1} (1972), 333--350.

\bibitem[Leh04]{Lehner:2002:Cumulants1}
Franz Lehner, \emph{Cumulants in noncommutative probability theory {I}.
  {N}oncommutative exchangeability systems}, Math. Zeitschr. \textbf{248}
  (2004), no.~1, 67--100, arXiv:math.CO/0210442.

\bibitem[NS97]{NicaSpeicher:1997:fourier}
Alexandru Nica and Roland Speicher, \emph{A ``{F}ourier transform'' for
  multiplicative functions on non-crossing partitions}, J. Algebraic Combin.
  \textbf{6} (1997), 141--160. \MR{98i:46070}

\bibitem[Pet90]{Petz:1990:definetti}
D{\'e}nes Petz, \emph{A de {F}inetti-type theorem with $m$-dependent states},
  Probab. Theory Related Fields \textbf{85} (1990), 65--72. \MR{91e:46089}

\bibitem[Pis00]{Pisier:2000:inequality}
Gilles Pisier, \emph{An inequality for $p$-orthogonal sums in non-commutative
  ${L}\sb p$}, Illinois J. Math. \textbf{44} (2000), 901--923. \MR{1 804 311}

\bibitem[Spe94]{Speicher:1994:multiplicative}
Roland Speicher, \emph{Multiplicative functions on the lattice of noncrossing
  partitions and free convolution}, Math. Ann. \textbf{298} (1994), 611--628.

\bibitem[Spe97]{Speicher:1997:universal}
\bysame, \emph{On universal products}, Free probability theory (Waterloo, ON,
  1995), Amer. Math. Soc., Providence, RI, 1997, pp.~257--266.

\bibitem[Spe98]{Speicher:1998:combinatorial}
\bysame, \emph{Combinatorial theory of the free product with amalgamation and
  operator-valued free probability theory}, Mem. Amer. Math. Soc. \textbf{132}
  (1998), no.~627, x+88.

\bibitem[Sta86]{Stanley:1986:Enumerative1}
Richard~P. Stanley, \emph{Enumerative combinatorics. {V}ol. {I}}, The Wadsworth
  \& Brooks/Cole Mathematics Series, Wadsworth \& Brooks/Cole Advanced Books \&
  Software, Monterey, CA, 1986, With a foreword by Gian-Carlo Rota.
  \MR{87j:05003}

\bibitem[Sta99]{Stanley:1999:Enumerative2}
\bysame, \emph{Enumerative combinatorics. {V}ol. 2}, Cambridge Studies in
  Advanced Mathematics, vol.~62, Cambridge University Press, Cambridge, 1999.
  \MR{2000k:05026}

\bibitem[St{\o}69]{Stormer:1969:symmetric}
Erling St{\o}rmer, \emph{Symmetric states of infinite tensor products of
  {$C\sp{\ast} $}-algebras}, J. Functional Analysis \textbf{3} (1969), 48--68.
  \MR{39 \#3327}

\bibitem[SW97]{SpeicherWoroudi:1997:boolean}
Roland Speicher and Reza Woroudi, \emph{Boolean convolution}, Free probability
  theory (Waterloo, ON, 1995), Amer. Math. Soc., Providence, RI, 1997,
  pp.~267--279.

\bibitem[Voi95]{Voiculescu:1995:operations}
Dan Voiculescu, \emph{Operations on certain non-commutative operator-valued
  random variables}, Ast{\'e}risque (1995), no.~232, 243--275, Recent advances
  in operator algebras (Orl{\'e}ans, 1992).

\end{thebibliography}

\providecommand{\bysame}{\leavevmode\hbox to3em{\hrulefill}\thinspace}
\providecommand{\MR}{\relax\ifhmode\unskip\space\fi MR }
\providecommand{\MRhref}[2]{%
  \href{http://www.ams.org/mathscinet-getitem?mr=#1}{#2}
}
\providecommand{\href}[2]{#2}

\end{document}